\def\date{25.9.06}

\documentclass[12pt]{article}

\usepackage{amssymb}
\usepackage{amsmath}

\font\tengoth=eufm10 at 10pt
\font\sevengoth=eufm7 at 6pt
\newfam\gothfam
\textfont\gothfam=\tengoth
\scriptfont\gothfam=\sevengoth

\def\frak{\fam\gothfam\tengoth}

\def\fL{{\frak L}}

\def\fa{{\frak a}}
\def\fb{{\frak b}}

\def\fe{{\frak e}}

\def\fg{{\frak g}}

\def\fk{{\frak k}}
\def\fl{{\frak l}}
\def\fm{{\frak m}}
\def\fn{{\frak n}}

\def\fq{{\frak q}}
\def\fp{{\frak p}}
\def\fr{{\frak r}}
\def\fs{{\frak s}}

\renewcommand{\:}{\colon}
\newcommand{\1}{{\bf 1}}

\newcommand{\g}{{\mathfrak g}}
\newcommand{\s}{{\mathfrak s}}
\newcommand{\h}{{\mathfrak h}}
\newcommand{\z}{{\mathfrak z}}

\renewcommand{\phi}{\varphi}

\newcommand{\subeq}{\subseteq}

\newcommand{\into}{\hookrightarrow}
\newcommand{\eps}{\varepsilon}

\newcommand{\N}{{\mathbb N}}
\newcommand{\Z}{{\mathbb Z}}
\newcommand{\R}{{\mathbb R}}
\newcommand{\C}{{\mathbb C}}
\newcommand{\Cns}{{\mathbb C}}

\newcommand{\Q}{{\mathbb Q}}

\renewcommand{\hat}{\widehat}
\renewcommand{\tilde}{\widetilde}

\renewcommand{\L}{\mathop{\bf L{}}\nolimits}

\newcommand{\GL}{\mathop{{\rm GL}}\nolimits}

\newcommand{\SL}{\mathop{{\rm SL}}\nolimits}

\newcommand{\gl}{\mathop{{\fg\fl}}\nolimits}

\newcommand{\lsl}{\mathop{{\fs\fl}}\nolimits}

\newcommand{\ad}{\mathop{{\rm ad}}\nolimits}
\newcommand{\Ad}{\mathop{{\rm Ad}}\nolimits}

\newcommand{\Hom}{\mathop{{\rm Hom}}\nolimits}

\newcommand{\Aut}{\mathop{{\rm Aut}}\nolimits}

\newcommand{\id}{\mathop{{\rm id}}\nolimits}

\newcommand{\rad}{\mathop{{\rm rad}}\nolimits}
\renewcommand{\dim}{\mathop{{\rm dim}}\nolimits}
\newcommand{\im}{\mathop{{\rm im}}\nolimits}

\newcommand{\Rarrow}{\Rightarrow}
 
\newcommand{\oline}{\overline}

\newcommand{\res}{\vert}

\newcommand{\sL}{{\lsl}}

\newcommand{\Spec}{{\rm Spec}}

\newcommand{\ssssarr}{\hbox to 15pt{\rightarrowfill}}
\newcommand{\sssarr}{\hbox to 20pt{\rightarrowfill}}
\newcommand{\ssarr}{\hbox to 30pt{\rightarrowfill}}
\newcommand{\sarr}{\hbox to 40pt{\rightarrowfill}}
\newcommand{\arr}{\hbox to 60pt{\rightarrowfill}}
\newcommand{\larr}{\hbox to 60pt{\leftarrowfill}}
\newcommand{\Arr}{\hbox to 80pt{\rightarrowfill}}

\newcommand{\sssmapright}[1]{\smash{\mathop{\sssarr}\limits^{#1}}}

\def\theoremname{Theorem}
\def\propositionname{Proposition}
\def\corollaryname{Corollary}
\def\lemmaname{Lemma}
\def\remarkname{Remark}
\def\conjecturename{Conjecture} 

\def\definitionname{Definition}
\def\exercisename{Exercise}
\def\examplename{Example}
\def\examplesname{Examples}
\def\problemname{Problem}
\def\problemsname{Problems}
\def\proofname{Proof}

\def\satzname{Satz} 
\def\koroname{Korollar}
\def\folgname{Folgerung}
\def\bemerkname{Bemerkung}
\def\aufgname{Aufgabe}

\def\beisname{Beispiel}
\def\beissname{Beispiele}
\def\bewname{Beweis}

\def\@thmcounter#1{\noexpand\arabic{#1}}
\def\@thmcountersep{}
\def\@begintheorem#1#2{\it \trivlist \item[\hskip 
\labelsep{\bf #1\ #2.\quad}]}
\def\@opargbegintheorem#1#2#3{\it \trivlist
      \item[\hskip \labelsep{\bf #1\ #2.\quad{\rm #3}}]}
\makeatother
\newtheorem{theo}{\theoremname}[section]
\newtheorem{propo}[theo]{\propositionname}
\newtheorem{coro}[theo]{\corollaryname}
\newtheorem{lemm}[theo]{\lemmaname}

\newenvironment{theorem}{\begin{theo}\it}{\end{theo}}

\newenvironment{proposition}{\begin{propo}\it}{\end{propo}}

\newenvironment{corollary}{\begin{coro}\it}{\end{coro}}

\newenvironment{lemma}{\begin{lemm}\it}{\end{lemm}}

\newtheorem{rem}[theo]{\remarkname}

\newenvironment{remark}{\begin{rem}\rm}{\end{rem}}

\newtheorem{stepnow}[theo]{}

\newtheorem{defin}[theo]{\definitionname} 

\newenvironment{definition}{\begin{defin}\rm}{\end{defin}}

\newtheorem{exer}[theo]{\exercisename}

\newtheorem{ex}[theo]{\examplename}

\newenvironment{example}{\begin{ex}\rm}{\end{ex}}

\newtheorem{exs}[theo]{\examplesname}

\newtheorem{conj}[theo]{\conjecturename}

\newtheorem{pr}[theo]{\problemname}

\newenvironment{problem}{\begin{pr}\rm}{\end{pr}}

\newtheorem{prs}[theo]{\problemsname}

\newcommand{\qed}{{\unskip\nobreak\hfil\penalty50\hskip .001pt \hbox{}
          \nobreak\hfil
          \vrule height 1.2ex width 1.1ex depth -.1ex
           \parfillskip=0pt\finalhyphendemerits=0\medbreak}\rm}

{\hfill\qed\end{trivlist}}
\newenvironment{proof}{\begin{trivlist}\item[\hskip%
\labelsep{\bf\proofname.\quad}]}%
{\hfill\qed\end{trivlist}}

\newenvironment{Proof*}{\begin{trivlist}\item[\hskip%
\labelsep{\bf\proofname.\quad}]}%
{\end{trivlist}}

{\hfill\qed\end{trivlist}}

\newenvironment{beweis*}{\begin{trivlist}\item[\hskip%
\labelsep{\bf\bewname.\quad}]}%
{\end{trivlist}}

\newtheorem{satzn}[theo]{\satzname}

\newtheorem{koro}[theo]{\koroname}

\newtheorem{folg}[theo]{\folgname}

\newtheorem{bem}[theo]{\bemerkname}

\newtheorem{aufg}[theo]{\aufgname}

\newtheorem{aufgn}[theo]{\aufgname}

\newtheorem{beis}[theo]{\beisname}

\newtheorem{beiss}[theo]{\beissname}

\begin{document} 

\newcommand{\tor}{\mathop{{\rm tor}}\nolimits}
\newcommand{\Spann}{\mathop{{\rm span}}\nolimits}
\newcommand{\comp}{\mathop{{\rm comp}}\nolimits}
\newcommand{\evol}{\mathop{{\rm evol}}\nolimits}
\newcommand{\der}{\mathop{{\rm der}}\nolimits}
\newcommand{\dist}{\mathop{{\rm dist}}\nolimits}
\newcommand{\proj}{\mathop{{\rm pr}}\nolimits}
\newcommand{\mot}{\mathop{{\mathfrak{mot}}}\nolimits}
\newcommand{\osc}{\mathop{{\mathfrak{osc}}}\nolimits}

\newcommand{\dl}{{\displaystyle \lim_{\longrightarrow}}\ }
\newcommand{\pl}{{\displaystyle\lim_{\longleftarrow}}\ }
\newcommand{\indlim}{\dl} 
\newcommand{\prolim}{\pl} 
\newcommand{\trile}{\trianglelefteq}

\renewcommand{\R}{{\mathbb R}}
\renewcommand{\C}{{\mathbb C}}
\renewcommand{\fL}{{\mathfrak L}}



\def\alert#1{\hfill\break\vbox{\smallskip
             \noindent\vrule\vbox{\hsize=.99\hsize
             \noindent\hrule\par
             \noindent\vskip3pt\hbox{\hskip3pt
             \vbox{\hsize=.98\hsize\noindent#1}\hskip3pt}\vskip3pt\par
             \noindent\hrule}\vrule}}

\def\alertshort#1{\hfill\break\vbox{\smallskip
             \noindent\vrule\vbox{\hsize=.75\hsize
             \noindent\hrule\par
                    \noindent\vskip3pt\hbox{\hskip3pt
             \vbox{\hsize=.98\hsize\noindent#1}\hskip3pt}\vskip3pt\par
             \noindent\hrule}\vrule}}



\title{Pro-Lie groups which are infinite-dimensional Lie groups} 

\author{K.\ H.\ Hofmann and K.-H. Neeb}

\maketitle


\begin{abstract} A pro-Lie group is a projective limit of a family of 
finite-dimensional Lie groups. In this note we show that 
a pro-Lie group $G$ is a Lie group in the sense that its topology 
is compatible with a smooth manifold structure for which the group 
operations are smooth if and only if $G$ is locally contractible. We also 
characterize the corresponding pro-Lie algebras in various ways. 
Furthermore, we characterize those pro-Lie groups 
which are locally exponential, 
that is,  they are Lie groups with a smooth exponential function which 
maps a zero neighborhood in the Lie algebra 
diffeomorphically  onto an open identity neighborhood of the group.
\end{abstract} 

\section{Introduction.}

There are several natural ways to extend the theory of finite-dimensional Lie groups 
to larger classes of groups. From a topological perspective, the closest relatives 
to finite-dimensional Lie groups are locally compact groups. According to a classical theorem 
of Yamabe, each locally 
compact group contains an open subgroup which is a projective limit of finite-dimensional 
Lie groups; let us call these groups {\it pro-Lie groups}. Hence the local structure 
of a locally compact group is that of a pro-Lie group, but, as examples such as 
the topological product $\R^\N$ show, not every pro-Lie group is locally compact. 
In the recent monograph \cite{HoMo06} Hofmann and Morris develop  
an effective Lie theory for 
the whole class of pro-Lie groups, which contains in particular the 
Lie theory of locally compact groups. 
Throughout this theory, pro-Lie groups are considered as topological groups with 
additional structural features, but not as Lie groups in a differentiable setting. 
However, there is a Lie functor $\fL$ assigning to each pro-Lie group a pro-Lie algebra, 
that is,  a projective limit of finite-dimensional Lie algebras. 
Projective limits of finite-dimensional Lie groups occur 
quite naturally under various aspects in the literature, see f.i., 
Lewis \cite{Lew39}, Kuranishi \cite{Kur59}, Sternberg \cite{Ste61} and Omori \cite{Omo80}. 

On the other hand, quite generally, we say that 
a {\it Lie group} is a group $G$, endowed with the structure of a 
manifold modelled on a locally convex space, such that the group operations 
on $G$ are smooth. We write $\1$ for the identity element of $G$. 
For any Lie group $G$, the tangent space $T_\1(G)$ can be identified with 
the space of left invariant vector fields, hence inherits the structure 
of a {\it locally convex Lie algebra}, that is,  a locally convex space with a 
continuous Lie bracket. We write $\L(G) := (T_\1(G), [\cdot,\cdot])$. 
A smooth map $\exp_G \: \L(G) \to G$ is said to be an {\it exponential 
function} if for each $x \in \L(G)$, the curve $\gamma_x(t) := \exp_G(tx)$ is a 
homomorphism $\R \to G$ with $\gamma_x'(0) = x$. Presently, all known 
Lie groups modelled on complete locally convex spaces possess an exponential 
function. For Banach--Lie groups, its existence follows from the 
theory of ordinary differential equations in Banach spaces. 
A Lie group $G$ is called {\it locally exponential}, if 
it has an exponential function mapping an open $0$-neighborhood in 
$\L(G)$ diffeomorphically onto an open neighborhood of $\1$ in $G$. 
For more details, we refer to Milnor's lecture notes \cite{Mil84}, 
the second author's recent survey \cite{Ne06} or his extensive monograph with Gl\"ockner
\cite{GN06}. 

\bigskip

It is the goal of the present paper to describe how the theory of pro-Lie groups 
intersects the theory of 
Lie groups in this sense. Clearly, any Lie group which is locally compact, 
is modelled on a finite-dimensional space, hence is a finite-dimensional Lie group. 


There are three 
natural questions to be answered for a pro-Lie group $G$: 
\begin{description}
\item[\rm(Q1)] When does $G$ carry a Lie group structure compatible with its topology? 
\item[\rm(Q2)] To which extent is the Lie group structure on a pro-Lie group unique? 
\item[\rm(Q3)] Suppose that the pro-Lie group 
$G$ carries a Lie group structure, when is $G$ locally 
exponential, 
that is, when is its exponential function a local diffeomorphism in~$0$?
\end{description}

In the process of presenting our answers, the concept of a space to be {\it locally contractible} plays an important role. We say that a topological group $G$ is locally
contractible if {\bf1} has arbitrarily small neighborhoods which are contractible
in $G$. Since several different
definitions are possible and feasible, we collected precise definitions in Appendix 10
below. Given this concept, we answer question (Q1) by characterizing Lie groups
among pro-Lie groups
in a purely topological fashion,
as follows (Theorems~\ref{thm:10} and \ref{thm:3.2}). 


\begin{theorem} \label{thm1} A pro-Lie group $G$ carries a Lie group 
structure compatible with its topology if and only if it is locally contractible. 
\end{theorem}

However, a more accessible characterization in terms of 
the Lie algebra of $G$ would be highly desirable.
Unfortunately, 
there might be several pro-Lie groups $G$ with isomorphic Lie algebras, some of 
which might be 
locally contractible while 
others are not. But if 
a pro-Lie group 
$G$ is locally 
contractible, then it has a universal covering group 
which still is pro-Lie  and 
locally contractible. 
We call a topological space $1$-connected if it is arcwise 
connected and has a trivial fundamental group. 
A $1$-connected pro-Lie group is completely determined by 
its Lie algebra, so that 
a Lie algebraic answer to question (Q1) 
 has to characterize those 
pro-Lie algebras $\g$ for which the 
universal $1$-connected group $\Gamma(\g)$ 
attached to it by Lie's Third Theorem for pro-Lie groups 
is locally contractible.  

\begin{theorem} \label{thm2} 
For a pro-Lie algebra $\g$, the following are equivalent: 
  \begin{description}
    \item[\rm(1)]  $\g$ is the Lie algebra of a locally convex Lie group $G$ 
      with smooth exponential function. 
    \item[\rm(2)] $\g$ has a Levi decomposition 
      $\g \cong \fr \rtimes \s$, where only finitely many factors in 
      $\s \cong \prod_{j \in J} \s_j$ are not isomorphic to $\sL_2(\R)$. 
    \item[\rm(3)] The corresponding $1$-connected 
universal group $\Gamma(\g)$ is locally contractible. 
    \item[\rm(4)] The maximal compact subgroups of $\Gamma(\g)$ are Lie groups. 
    \item[\rm(5)] There exists a locally contractible pro-Lie group $G$ with 
$\fL(G) \cong \g$. 
  \end{description}
\end{theorem}

This theorem follows from Theorem~\ref{thm:10} ($(3)\Leftrightarrow(4)$), 
Theorem~\ref{thm:2.5} ($(2)\Leftrightarrow(3)\Leftrightarrow(5)$) and 
Theorem~\ref{thm:nec-smooth}. 

In the following, we shall call pro-Lie algebras satisfying the equivalent conditions 
in the preceding theorem {\it smooth}. 

The uniqueness of Lie group structures can be treated in two essentially 
different ways. One is to shift the focus to the corresponding pro-Lie algebra $\g$ 
and ask for the uniqueness of corresponding $1$-connected Lie groups. 
Such uniqueness results are available under the assumption 
that the Lie group is 
regular (see Definition~\ref{def:reg} below). In this 
direction we shall show 

\begin{theorem} \label{thm3} For a pro-Lie algebra $\g$, the following assertions hold: 
  \begin{description}
  \item[\rm(1)] If $G$ is a Lie group with a smooth exponential function and 
$\L(G) = \g$, then $\g$ is smooth. 
  \item[\rm(2)] If $\g$ is smooth, then there exists a unique 
$1$-connected regular Lie group, which is 
isomorphic to $\Gamma(\g)$ as a topological group.  
 \item[\rm(3)] If $G$ is any connected regular Lie group for which $\g = \L(G)$, 
then $G$ is a quotient of $\Gamma(\g)$ by a discrete central subgroup. A subgroup 
$D \subeq Z(\Gamma(\g))$ is discrete if and only if it is finitely generated 
and its intersection with the identity component $Z(\Gamma(\g))_0 \cong \z(\g)$ 
is discrete. 
  \end{description}
\end{theorem}

Part (1) follows from Theorem~\ref{thm:nec-smooth}; Part (2) is a consequence of 
Proposition~\ref{prop:reglie} and Corollary~\ref{cor:reg-iso}, and Part (3) 
is taken from Theorems~\ref{thm:10} and \ref{thm:3.8}. 

Part (3) of the preceding theorem provides in particular a quite tractable 
description of all connected regular Lie groups whose Lie algebra is pro-Lie. 

A second strategy to address the uniqueness question is to use local exponentiality 
instead of regularity of the Lie group under consideration. Here a Lie group is called 
{\it locally exponential} if it has an exponential function which is a local diffeomorphism. 
Local exponentiality is well compatible with topological group structures because 
continuous morphisms of locally exponential Lie groups are automatically smooth, 
so that locally exponential Lie groups form a full subcategory of the 
category of topological groups (cf.\ \cite{GN06}). 
Thus it makes sense to call a topological 
group {\it locally exponential} if it carries a locally exponential Lie group 
structure compatible with the given topology. 

\begin{definition} \label{def:2.2} A locally convex Lie algebra $\g$ 
is called {\it locally exponential} if there exists a circular 
convex open $0$-neighborhood $U \subeq \g$ and a smooth map 
$$ U \times U \to \g, \quad (x,y) \mapsto x * y $$
satisfying:
\begin{description}
\item[\rm(E1)] $x * (y*z) = (x*y)*z$ if $x*y, y*z \in U$. 
\item[\rm(E2)] $x * 0 = 0 * x = x$. 
\item[\rm(E3)] For $x \in U$ and $|s|,|t|\leq 1$, we have $sx * tx = (s+t)x.$ 
\item[\rm(E4)] The second order Taylor polynomial of $*$ in $0$ is given by 
$x + y + \frac{1}{2}[x,y].$ 
\end{description}

A locally convex Lie algebra $\g$ is called {\it exponential} if 
the above conditions are satisfied for $U = \g$. 
In view of (E3), this means that $(\g,*)$ is a Lie group whose 
exponential function $\exp_\g$ coincides with $\id_\g$. 
\end{definition}

Since any local Lie group on an open subset of a locally convex space 
$V$ leads to a Lie algebra structure on $V$ (\cite{GN06}), 
condition (E4) only ensures that $\g$ is the Lie algebra of the 
corresponding local group. 

We have the following answer to question (Q3), which uses the concept of a locally 
exponential Lie algebra, defined in Definition~\ref{def:2.2} 
(Theorem~\ref{thm:locexp}, Corollary~\ref{cor:locexp}): 

\begin{theorem} \label{thm4} A pro-Lie group $G$ is locally exponential if and only if 
it is locally contractible and $\fL(G)$ is a locally exponential Lie algebra. 
A pro-Lie algebra $\g$ is locally exponential if and only if 
the set of $\exp$-regular points, that is,  the set of all $x \in \g$ for which 
$$ \Spec(\ad x) \cap 2\pi i \Z = \{0\}, $$
is a $0$-neighborhood.   
\end{theorem}

In the framework of things considered here, a proof of 
Theorem~\ref{thm4} is harder than one might surmise.

The structure of the paper is as follows. Section 2 surveys some key results 
on pro-Lie groups, mostly cited from \cite{HoMo06}, 
and in Section 3 we discuss 
their local contractibility. Section 4 contains various characterizations of 
smooth pro-Lie algebras, whereas Lie group structures on the corresponding groups 
are treated in Section~5. In Section 6 we recall some results on locally 
exponential Lie algebras from \cite{GN06} and in Section 7 we finally 
prove Theorem~\ref{thm4}. 

\section{Pro-Lie groups and their Lie algebras} 

In this section, we collect some of the key results of \cite{HoMo06} 
concerning pro-Lie groups. 
Clearly, arbitrary products of finite-dimensional Lie groups, such as 
$$ \R^J, \quad \Z^J, \quad \SL_2(\R)^J $$
for an arbitrary set $J$, are pro-Lie groups. 
The following theorem provides an abstract characterization of pro-Lie groups:

\begin{theorem} \label{thm:2.0} {\rm(\cite{HoMo06})} A 
topological group $G$ is a pro--Lie group if and only if 
it is isomorphic to a closed subgroup of a product of finite-dimensional Lie 
groups. In particular, closed subgroups of pro-Lie groups are pro-Lie groups. 
\end{theorem}

Since projective limits are defined as certain 
closed subgroups of pro\-ducts, one implication 
of the preceding theorem is trivial; the converse is the interesting part. 

\begin{theorem} \label{thm:2.0b} {\rm(\cite{HoMo06})} If 
$G$ is a pro-Lie group, then there exists a 
filter basis ${\cal N}$ of closed normal subgroups $N \trile G$ 
for which $G/N$ is a finite-dimensional Lie group and $\lim {\cal N} = \1$. 
\end{theorem}

For the equivalence of these various equivalent possible definitions of pro-Lie groups
see \cite{HoMo06}, Theorem~3.39.

The key tool to the Lie theory of pro-Lie groups is the observation that 
they ``have a Lie algebra'' in the following sense: 

\begin{definition} \label{def:2.1} (\cite{HoMo06}, Definition~2.11) 
Let $G$ be a topological group 
and $\fL(G) = \Hom_c(\R,G)$ the set of continuous one-parameter groups, 
endowed with the compact open topology. 
We define a scalar multiplication on $\fL(G)$ by 
\begin{eqnarray}
  \label{eq:2.1}
(\lambda\alpha)(t) := \alpha(\lambda t) \quad \mbox{ for } \quad 
\lambda \in \R, \alpha \in \Hom_c(\R,G). 
\end{eqnarray}
We say that $G$ is {\it a topological group with Lie algebra} if 
for $\alpha,\beta \in \fL(G)$ 
the limits 
\begin{eqnarray}
  \label{eq:2.2}
(\alpha + \beta)(t) := \lim_{n \to \infty} 
\Big(\alpha(\frac{t}{n})\beta(\frac{t}{n})\Big)^n
\end{eqnarray}
and 
\begin{eqnarray}
  \label{eq:2.3}
[\alpha,\beta](t^2) := \lim_{n \to \infty} \Big(
\alpha(\frac{t}{n})\beta(\frac{t}{n})
\alpha(-\frac{t}{n})\beta(-\frac{t}{n})\Big)^{n^2}. 
\end{eqnarray}
exist in the compact open topology, define elements of ${\frak L}(G)$, 
addition and bracket are continuous maps 
${\frak L}(G) \times {\frak L}(G) \to {\frak L}(G),$
and with respect to the scalar multiplication (\ref{eq:2.1}), the addition (\ref{eq:2.2}), 
and the bracket (\ref{eq:2.3}),  
${\frak L}(G)$ is a real Lie algebra. This implies that $\fL(G)$ is a topological Lie algebra. 

For any topological group $G$, we define the {\it exponential function of $G$} by 
$$ \exp_G \: \fL(G) \to G, \quad \alpha \mapsto \alpha(1). $$
\end{definition} 

A crucial observation is that the class of topological groups with 
Lie algebra is closed under projective limits and that 
$$ \fL(\prolim G_j) \cong \prolim \fL(G_j), $$
as topological Lie algebras (cf.\ \cite{HoMo06}, Theorem\ 2.25). 

Let us call topological vector spaces of the form $\R^J$, $J$ a set, 
{\it weakly complete}. These are the dual spaces of the vector spaces 
$\R^{(J)}$, endowed with the weak-$*$-topology. This provides a duality between 
real vector spaces and weakly complete locally convex spaces, which implies 
in particular that each closed subspace of a weakly complete space 
is weakly complete and complemented. 
For a systematic treatment see \cite{HoMo06}, App.\ 2. 
 In particular, weakly complete spaces are 
nothing but the projective limits of finite-dimensional vector spaces. 
These considerations lead to: 

\begin{theorem} \label{thm:2.1} {\rm(\cite{HoMo06},
Corollary~4.21, 4.22)} 
Every pro-Lie group $G$ has a Lie algebra $\fL(G)$ which is a 
a projective limit of finite-dimensional Lie algebras, hence a weakly complete 
topological Lie algebra. The image of the exponential function 
ge\-ne\-rates a dense subgroup of the identity component $G_0$. 
\end{theorem}

In the following, we call projective limits of finite-dimensional 
Lie algebras {\it pro-Lie algebras}. 

In view of Theorem~\ref{thm:2.0}, the category of pro-Lie groups is closed under products 
and projective limits. These remarkable closedness properties lead 
to the existence of an adjoint functor $\Gamma$ for the Lie functor~$\fL$: 

\begin{theorem} \label{thm:2.2} {\rm(Lie's Third Theorem for Pro-Lie Groups; \cite{HoMo06}, Theorem\ 2.26)} 
The Lie functor $\fL$ from the category of pro-Lie groups 
to the category of pro-Lie algebras has a left adjoint $\Gamma$. It associates 
with each pro-Lie algebra $\g$ a $1$-connected pro-Lie group $\Gamma(\g)$ 
and a natural isomorphism $\eta_\g \: \g \to \fL(\Gamma(\g))$, such that for 
every morphism $\phi \: \g \to \fL(G)$ of pro-Lie algebras, $G$ a pro-Lie group, 
there exists a unique morphism 
$\phi' \: \Gamma(\g) \to G$ with 
$\fL(\phi') \circ \eta_\g = \phi$. 
\end{theorem}

\begin{remark} \label{rem:2.6} If $\g$ is finite-dimensional, then $\Gamma(\g)$ is a $1$-connected 
Lie group with Lie algebra $\g$. 

One can also show that $\Gamma$ preserves semidirect products 
(\cite{HoMo98}, \break Theorem~6.11). Moreover, if 
$\g = \prolim \g_j$ is a general projective limit, we also have 
$\Gamma(\g) \cong \prolim \Gamma(\g_j),$
which often provides an explicit 
description of $\Gamma(\g)$ in many cases (cf.\ \cite{HoMo06},
Chapters\ 6 and 8). 
\end{remark}

It is quite remarkable that the category of pro-Lie algebras 
permits us to 
develop a structure theory which is almost as strong as in finite dimensions. 
In particular, there is a Levi decomposition. To describe it, we call a 
pro-Lie algebra $\g$ {\it prosolvable} if it is a projective limit 
of finite-dimensional solvable Lie algebras: 

\begin{theorem}  \label{thm:2.4} 
{\rm(Levi decomposition; \cite{HoMo06}, Theorems 7.52 and 7.77)} 
Each pro-Lie algebra $\g$ contains a unique maximal prosolvable ideal 
$\fr = \rad(\g)$ and $\s := \g/\fr$ is a product $\prod_{j \in J} \s_j$ 
of finite-dimensional simple Lie algebras~$\s_j$. We further have a Levi splitting 
$\g \cong \fr \rtimes \s,$ 
and two Levi factors $\s_1$ and $\s_2$ are conjugate 
under an inner automorphism of the form $e^{\ad x}$.
\end{theorem}

\section{Locally contractible pro-Lie groups} 

Since local contractibility is clearly necessary for a topological group to 
carry a compatible Lie group structure, we devote the present section to the topological 
structure of the locally contractible pro-Lie groups, the main result 
being Theorem~\ref{thm:10}, saying that all these groups are homeomorphic to 
products of vector spaces and compact Lie groups. 

In \cite{HoMo06}, the following result on connected pro-Lie groups
is established in 12.81 and 12.82:

\begin{theorem} \label{thm:4} 
Each connected pro-Lie group $G$ is homeomorphic to a product
of a compact connected semisimple subgroup $S$ of $G$, 
a compact connected abelian subgroup $A$ of $G$,
and a family  of copies of $\R$.
Moreover, $A$ is in the normalizer of $S$ and $SA$ 
is a maximal compact subgroup of $G$.
\end{theorem}

It is also shown that  and 
every compact subgroup of $G$ has a conjugate contained in $SA$.
Incidentally, this shows, among other things, that
 each connected pro-Lie  group is homeomorphic to a product
of a compact connected semisimple group and a connected abelian
pro-Lie group. (For the structure of connected abelian pro-Lie groups
see \cite{HoMo06}, Chapter 5.)                          

If $G$ is locally contractible, then by Lemma~\ref{lem:a.3} (in the appendix), both $S$ and $A$
are locally contractible. This causes us to discuss locally contractible
compact connected semisimple groups and locally contractible
compact connected abelian groups.

We summarize what is known on the contractibility of homogeneous spaces
of compact groups. The following theorem of A. Borel and its proof
is published in the appendix of \cite{HoMs66}, pp.\ 306-310, notably
Theorem~4.3, p.~310. Recall that a space $X$
is called acyclic over a ring $R$ with identity, 
if the \v Cech cohomology ring
$H^*(X,R)$ is that of a singleton space, that is, $H^0(Z,R)\cong R$, and 
$H^n(X,R)=\{0\}$ for all $n>0$. Note that $H^0(X,R)=C(X,R)$, where
$R$ is viewed with the discrete topology, whence $X$ is connected
if and only if $H^0(X,R)\cong R$.

\begin{theorem} \label{thm:5} 
For a compact group $G$ and a closed subgroup
$H$, the following statements are equivalent:
  \begin{description}
  \item[\rm(i)] $G/H$ is singleton.
  \item[\rm(ii)] $G/H$ is acyclic over $\Q$ and over $\Z/2\Z$.
  \end{description}
\end{theorem}

This applies in particular to the case $H=\{\1\}$ and characterizes
the degeneracy of a compact group in cohomological terms.
As a corollary, we get (\cite{HoMo98}, p.~310, 4.4)

\begin{corollary} \label{cor:homspace}
A homogeneous space of a compact group is contractible
if and only if it is singleton. 
\end{corollary}

In particular, this applies to compact groups.

Now we address first the semisimple case:
\begin{proposition} \label{prop:semisim}
Let $S$ be a connected  compact semisimple
group. Then the following statements are equivalent:
\begin{description}
\item[\rm(i)] $S$ is a Lie group.
\item[\rm(ii)] $S$ is locally contractible. 
\end{description}
\end{proposition}

\begin{proof}
By Remark~\ref{rem:a.2} (in the appendix), that implication (i) implies (ii) 
is trivial and we have to prove that (ii) implies (i). 
So assume there is a neighborhood $U$ of the identity which is contractible
in $S$. Since every compact group is a projective limit of Lie
groups, there is a compact normal subgroup $N\subseteq U$ such that
$S/N$ is a Lie group.  The structure theory of compact connected 
semisimple groups in \cite{HoMo98}, Theorem~9.19(i) and (ii) allows us to derive 
that $N/N_0$ is finite and that $S/N_0$ is a Lie group. Thus we may assume
that $N$ is connected. From loc.\ cit.\ we deduce the existence of a compact normal
connected semisimple subgroup $L$ such that $S=NL$ and $N\cap L$ is
central finite. The morphism $\mu\colon N\times L\to S$, $\mu(n,x)=nx$
is a covering morphism. As such, it has the homotopy lifting property.

Since $N\subseteq U$ and $U$ is contractible in $S$, we have a
a contraction \break $H\colon N\times[0,1]\to S$ of $N$  in $S$, such that
$H(n,0)=n$, $H(n,1)=\1$ for all $n\in N$.
 By the homotopy lifting
 property, there is a $\tilde H\colon N\times [0,1]\to N\times L$
such that $\tilde H(n,0)=(n,\1)$, $\tilde H(n,1)=(\1,\1)$. That is, 
$N\times \{\1\}$ is contractible in $N\times L$. 
Let $j\colon N\to N\times L$ be the inclusion $n\mapsto (n,1)$. Then
in the homotopy category, $[j]=[\1]$, the homotopy class of the constant function
with value $\1$. Now let $\proj_N\colon N\times L\to N$ be the projection onto $N$. It 
satisfies $\proj_N\circ j
=\id_N,$ and thus $[\id_N]=[\1]$, that is, $N$ is contractible.

Now Corollary~\ref{cor:homspace} 
applies and shows that $N=\{\1\}$, whence $S\cong S/N$ is a Lie group.
\end{proof}

The next step is the abelian case.

\begin{proposition} \label{prop:abel} Let $A$ be a connected  compact abelian
group. Then the following statements are equivalent:
\begin{description}
\item[\rm(i)] $A$ is a Lie group.
\item[\rm(ii)] $A$ is locally contractible. 
\end{description}
\end{proposition}

\begin{proof} Again we only have to prove that (ii) implies (i).

Let $U$ be a neighborhood  of the identity which is contractible
in $A$.
By the Iwasawa local decomposition theorem (\cite{Iwa49}, see also \cite{HoMo06},
Corollary~13.20),  there is a closed subgroup
$N$ contained in $U$ and a continuous map $f\colon\R^n\to A$
such that $p\colon N\times\R^n\to A$, $p(n,v)=nf(v)$ is a
covering morphism. Since $N$ is contractible in $A$, we may argue
just as in the proof of Proposition~\ref{prop:semisim} that $N\times\{0\}$ is
contractible in $N\times\R^n$ , which in turn implies that $N$
is contractible. Thus $N$ is singleton by Corollary~\ref{cor:homspace}. Thus
we have a covering morphism $\R^n\to A$ which shows that $A$
is a Lie group, as asserted.
\end{proof}

\begin{corollary} \label{cor:9}
Any locally contractible connected compact group is a Lie group.
\end{corollary}

\begin{proof} By the Borel--Scheerer--Hofmann Splitting Theorem (\cite{HoMo98}, Theorem~9.39),
$G= G'\rtimes A$ with a closed abelian subgroup $A\cong G/G'$. If $G$ is locally contractible
so are $G'$ and $A$ by Lemma 10.3 in Appendix 10. Then Propositions~\ref{prop:semisim} 
and \ref{prop:abel} 
imply that $G'$ and $A$ are Lie groups, and so $G$ is a Lie group as asserted.
 \end{proof}

Now we have the ingredients to prove

\begin{theorem} \label{thm:10} 
For a connected pro-Lie group $G$ the following are equivalent: 
\begin{description}
\item[{\rm(1)}] $G$ is locally contractible. 
\item[{\rm(2)}] A maximal compact subgroup $C$ of $G$ is a Lie group. 
\item[{\rm(3)}] $\Gamma(\fL(G))$ is locally contractible and 
$G \cong \Gamma(\fL(G))/D$ for some discrete central subgroup~$D$. 
\end{description}

If these conditions are satisfied, then $G$ is homeomorphic to 
$\R^J\times C$ for a set $J$.
\end{theorem}

Recall that, if a locally compact topological group is an (a priori infinite-dimensional) 
Lie group, then it is modelled on a locally compact, hence finite-dimensional space and 
therefore it is finite-dimensional. 

\begin{proof} By Theorem~\ref{thm:4}, $G$ is homeomorphic to $\R^J\times C$,  
where  $C$ is a maximal compact subgroup of $G$, which proves the last statement. 

(1) $\Rarrow$ (2): By Lemma~\ref{lem:a.3} (in the appendix), $C$ is locally contractible,
so that Corollary~\ref{cor:9} implies that $C$ is a Lie group. 

(2) $\Rarrow$ (1) is clear. 

(3) $\Rarrow$ (2): By assumption, $q \: \Gamma(\fL(G)) \to G$ is a covering 
map, 
so that the local contractibility of $G$ directly follows from the local contractibility 
of $\Gamma(\fL(G))$. 

(2) $\Rarrow$ (3): Let $G$ be a locally contractible pro-Lie group 
with Lie algebra $\fL(G) = \g$. 
By Theorem~\ref{thm:10}, $G$ is locally connected. Therefore
it has a universal covering group $\tilde G$ and 
$G \cong \tilde G/D$ holds for a discrete central subgroup~$D$. By
\cite{HoMo06}, Theorem~8.21, we have $\tilde G\cong\Gamma(\g)$. 
In particular, $\tilde G$ is a pro-Lie group as well. Let $U$ be a connected identity
neighborhood of $G$ and $H\colon[0,1]\times U\to G$ a contraction of
$U$ in $G$. If $U$ is sufficiently small, then $\tilde G$
has an identity neighborhood $\tilde U$ which is mapped homeomorphically
onto $U$ under the universal covering morphism. Then $H$
lifts to a contraction $\tilde H\colon[0,1]\times\tilde U\to\tilde G$ of
$\tilde U$ in $\tilde G$. Hence $\Gamma(\g)\cong\tilde G$ is locally 
contractible, and so (3) is proved.
\end{proof}

\section{Smooth pro-Lie algebras} 

The goal of this section are various characterizations of smooth pro-Lie algebras. 

\begin{definition} \label{def:smooth} 
We call a pro-Lie algebra $\g$ {\it smooth} if the corresponding 
simply connected universal pro-Lie group $\Gamma(\g)$ is locally contractible. 
\end{definition}

We shall see in the following section that smoothness of a pro-Lie algebra 
is equivalent to the existence of a Lie group structure on a corresponding pro-Lie 
group; justifying the terminology. 

\begin{definition} \label{def:reg} A subalgebra $\fk$ of a finite-dimensional real Lie algebra 
$\g$ is said to be {\it compactly embedded} if $\exp(\ad \fk)$ is contained 
in a compact subgroup of $\Aut(\g)$. 

To formulate the corresponding concept for pro-Lie algebras, we call a
 finite-dimensional module $(V, \rho_V)$ of a Lie algebra $\fk$ {\it compact} if 
$\exp(\rho_V(\fk))$ is contained in a compact subgroup of $\GL(V)$. 
A projective limit of finite-dimensional compact modules is called 
{\it pro-compact}. Now, a subalgebra $\fk$ of a pro-Lie algebra $\g$ 
is said to be {\it compactly embedded} if $\g$ is a pro-compact $\fk$-module. 
\end{definition}

It is shown in \cite{HoMo06} (Theorems~12.15 and 12.27) that each compactly 
embedded subalgebra is contained in a maximal one and that two maximal 
ones are conjugate under (inner) automorphisms of $\g$. 

The following lemma belongs to the folklore of finite-dimensional Lie theory; 
we recall its proof for the sake of completeness and later applications below. 

\begin{lemma} \label{lem:findim-cont} Let 
$G$ be a $1$-connected finite-dimensional Lie group
and $\fk \subeq \g := \L(G)$ a maximal compactly embedded subalgebra. 
Then the following are equivalent: 
\begin{description}
\item[\rm(1)] $G$ is contractible. 
\item[\rm(2)] All compact subgroups of $G$ are trivial. 
\item[\rm(3)] All simple ideals in $\g/\rad(\g)$ are isomorphic to $\sL_2(\R)$. 
\item[\rm(4)] $\fk$ is abelian. 
\end{description}
\end{lemma}

\begin{proof} We choose a Levi decomposition 
$\g = \fr \rtimes \s$ for which 
$\fk = \fk_\fr \oplus \fk_\s$ holds for 
$\fk_\fr := \fk \cap \fr$ and $\fk_\fs := \fk \cap \fs$ and recall 
that $\fk_\s$ is maximal compactly embedded in $\s$ (\cite{Ne99}, Proposition~VII.1.9).
Let $\s =\s_1\oplus\cdots\oplus \s_n$ be the decomposition into simple ideals 
and observe that $\fk_\s$ is adapted to this decomposition in the sense that 
$\fk_\s =\fk_1\oplus\cdots\oplus \fk_n$ for $\fk_j := \fk_\fs \cap \fs_j$
and $\fk_j$ is a maximal compactly embedded subalgebra of $\s_j$.

As $G$ is 1-connected, we know that $G\cong R\rtimes(S_1\times\cdots\times S_n)$
with the 1-connected radical $R$ and 
a finite sequence of 1-connected simple Lie groups $S_i$.
As a 1-connected solvable Lie group, $R$ is contractible. Therefore (1) is equivalent
to
\item{(1$'$)} $S_j$ is contractible for each $j=1,\dots,n$.

Let $S_j=K_jA_jN_j$ be the Iwasawa decomposition, where $\L(K_j)\cong\fk_j$. Since
the factor $A_jN_j$ is always diffeomorphic to a euclidean space, it is contractible.
If follows that (1$'$) is equivalent to
\item{(1$''$)} $K_j$ is contractible for each $j=1,\ldots,n$.

Now $K_j$ is compactly embedded and thus is of the form $C_j\times \R^{n_j}$ with
a compact connected Lie group $C_j = K_j'$. Consequently Corollary~\ref{cor:homspace} shows
that $K_j$ is contractible if and only if $C_j=\{\1\}$. This says that the maximal compact
subgroups of $S_j$ are trivial for each $j$. 
Therefore (1) and (2) are equivalent, that is,  $G$ is contractible if and only if 
the maximal compact subgroup $C := C_1 \times \cdots \times C_n$ is trivial. 
In view of $C = K'$, $K$ is abelian if and only if $C$ is  trivial, so that 
(2) and (4) are equivalent. 

It remains to see that (3) is equivalent to (4). 
Assume (4), that is,  that $\fk$ is abelian. 
Let $\fp_j$ be the orthogonal complement of $\fk_j$ in $\s_j$ with respect to the Killing 
form of $\s_j$. 
Since $\fk_j$ is abelian, it is one-dimensional (cf.\ \cite{Hel78}, Theorem~6.2), so that 
the simple $\fk_j$-module $\fp_j$ is $2$-dimensional, which implies that 
$\dim \s_j = 3$ and hence that $\s_j \cong \sL_2(\R)$, because 
it is a non-compact simple Lie algebra. To see that, conversely, (3) implies (4), 
note that $\s_j \cong \sL_2(\R)$ implies $\fk_j \cong \R$ and hence that $\fk$ is abelian. 
This completes the proof.
\end{proof} 

\begin{remark} \label{rem:bch-mult} 
If $\fn = \prolim \fn_j$ is a pronilpotent Lie algebra, 
where all $\fn_j$ are nilpotent, then 
$\Gamma(\fn_j) \cong (\fn_j, *)$, where $*$ denotes the (polynomial) BCH multiplication 
on $\fn_j$. Since each connecting map $\phi_{ij} \: \fn_j \to \fn_i$ induces a 
homomorphism of groups $(\fn_j,*) \to (\fn_i,*)$, we obtain 
$\Gamma(\fn) \cong \prolim (\fn_j,*) \cong (\fn,*),$
where $*$ is given by the BCH-series, which converges on $\fn \times \fn$ because 
it converges on each finite-dimensional quotient. 
\end{remark} 

Recall from Definition~\ref{def:smooth} that a pro-Lie algebra $\g$ is called 
smooth if $\Gamma(\g)$ is locally contractible. 

\begin{theorem} \label{thm:2.5} For 
a pro-Lie algebra $\g$ with prosolvable radical $\fr$, the following are equivalent: 
  \begin{description}
  \item[\rm(1)] $\g$ is smooth. 
  \item[\rm(2)] The semisimple Lie algebra $\s := \g/\fr$ 
contains only finitely many simple ideals not isomorphic to $\sL_2(\R)$. 
  \item[\rm(3)] The $1$-connected pro-Lie group 
$\Gamma(\g)$ is a topological manifold modelled on~$\g$. 
  \item[\rm(4)] There exists a 
locally contractible pro-Lie group 
with Lie algebra $\g$. 
  \end{description}
\end{theorem}

\begin{proof} (1) $\Rarrow$ (2): In view of the Theorem~\ref{thm:2.4}, 
$\Gamma(\g) \cong \Gamma(\fr) \rtimes \Gamma(\s),$ where 
$\Gamma(\fr)$ is homeomorphic to $\fr$, hence 
contractible and $\Gamma(\s) \cong \prod_{j \in J} \Gamma(\s_j)$. 

Let $U$ be an identity neighborhood in $\Gamma(\g)$ which is contractible in $\Gamma(\g)$. 
Then there is a cofinite subset $J_1 \subeq J$ such that 
$\prod_{j \in J_1} \Gamma(\s_j) \subeq U.$
For each $j \in J$, we have morphisms of topological groups 
$\alpha_j \: \Gamma(\s_j) \to \Gamma(\g)$, $\beta_j \: \Gamma(\g) \to \Gamma(\s_j)$ with $\beta_j \circ \alpha_j = \id_{\Gamma(\s_j)}$ 
and for $j \in J_1$ the contractibility of $U$ in $\Gamma(\g)$ implies that 
$\beta_j$ is homotopic to a constant map, hence that 
$\id_{\Gamma(\s_j)}$ is homotopic to a constant map, 
that is,  $\Gamma(\s_j)$ is contractible. 
From Lemma~\ref{lem:findim-cont} we now derive that $\s_j \cong \sL_2(\R)$. 

(2) $\Rarrow$ (3): Let $J_0 := \{ j \in J \: \s_j \cong \sL_2(\R)\}$, which is a 
cofinite subset of~$J$. Then 
$$ \Gamma(\g) \cong \Gamma(\fr) \rtimes \Gamma(\s) \cong \big(\Gamma(\fr) 
\rtimes \Gamma(\s_0)\big) \rtimes \Gamma(\s_1), $$
where $\s_1 := \prod_{j \not\in J_0} \s_j$ is finite-dimensional and 
$\s_0 := \prod_{j \in J_0} \s_j \cong \sL_2(\R)^{J_0}$. 
Since $\tilde\SL_2(\R)$ is homeomorphic to $\R^3$, the subgroup 
$\Gamma(\s_0)$ is homeomorphic to $(\R^3)^{J_0}$. 

Let $\fe \subeq \fr$ be a closed vector space complement 
of $\fn := \oline{[\fr,\fr]}$ (\cite{HoMo06}, A.2.12(a)). 
The Lie algebra $\fn$ is pronilpotent  
because all images of this subalgebra in finite-dimensional 
quotients of $\fr$ are nilpotent. Moreover, the map 
\begin{eqnarray}
  \label{eq:3.1}
\Phi \: \fn \times \fe \mapsto \Gamma(\fr), \quad (x,y) \mapsto \exp_{\Gamma(\fr)}(x)
\exp_{\Gamma(\fr)}(y). 
\end{eqnarray}
is a homeomorphism (\cite{HoMo06}, Theorem~8.13). 
In particular, $\Gamma(\fr)$ is homeomorphic to $\fn \times \fe \cong \fr$. 

We conclude that $\Gamma(\fr)$, $\Gamma(\s_0)$ and $\Gamma(\s_1)$ 
are topological manifolds, which implies the assertion. 
  
(3) $\Rarrow$ (4) is trivial. 

(4) $\Rarrow$ (1) follows from Theorem~\ref{thm:10}. 
\end{proof}

\begin{remark} There is an alternative argument for the implication 
(1) $\Rarrow$ (2) in the preceding theorem, based on Theorem~\ref{thm:10}. 

The local contractibility of $\Gamma(\g)$ implies that its maximal compact 
subgroup is finite-dimensional. Therefore at most finitely many of the 
groups $\Gamma(\s_j)$ contain non-trivial compact subgroups, which is equivalent 
to $\s_j \not\cong \sL_2(\R)$. 
\end{remark}

The philosophy that the maximal compact subgroups of a finite-dimensional 
Lie group determine its topological behavior carries over to pro-Lie groups.
In the context of the preceding theorem, it leads to a formulation of smoothness
in terms of maximal compactly embedded subalgebras. 

\begin{proposition} \label{prop:2.12} Let 
$\fk \leq \g$ be a maximal compactly embedded subalgebra. 
Then the following are equivalent: 
\begin{description}
\item[\rm(1)] $\g$ is smooth. 
\item[\rm(2)] $\fk$ is smooth. 
\item[\rm(3)] $\fk$ is nearly abelian, that is,  
its commutator subalgebra is finite-dimensional. 
\end{description}
\end{proposition}

\begin{proof} Let $\g = \fr \rtimes \s$ be a Levi decomposition and 
$\fk_\s \subeq \s$ a maximal compactly embedded subalgebra. 
We claim that $\fk_\s$ can be enlarged to a maximal compactly embedded subalgebra 
$\fk$ of $\g$ satisfying 
\begin{eqnarray}
  \label{eq:levi-comp}
\fk = \fk_\fr \oplus \fk_\s \quad \mbox{ with } \quad 
\fk_\fr \subeq \z(\fk). 
\end{eqnarray}
In fact, the finite-dimensional representation theory of semisimple Lie algebras 
(Weyl's trick) implies that $\g$ is a compact $\fk_\s$-module, that is,  $\fk_\s$ is compactly 
embedded in $\g$. Now \cite{HoMo06}, Theorem~12.15,   
implies that $\fk_\s$ is contained in a maximal compactly embedded subalgebra $\fk$ of $\g$. 
Then the projection $\g = \fr \rtimes \s \to \s$ maps $\fk$ into a compactly embedded subalgebra 
of $\s$, and the maximality of $\fk_\s$ thus shows that 
$\fk \subeq \fr \rtimes \fk_\s$, hence that $\fk$ is adapted to the Levi decomposition of $\g$. 
As $\fk_\fr$ is a prosolvable ideal of the pro-compact Lie algebra $\fk$, it is central. 
This proves our claim. 

Next we recall from \cite{HoMo06}, Theorem~12.27,  
that all maximal compactly embedded subalgebras 
are conjugate under inner automorphisms of $\g$.  
We may therefore assume that $\fk$ satisfies~(\ref{eq:levi-comp}). 

Let $\s = \prod_{j \in J} \s_j$ be the decomposition into simple ideals 
and observe that $\fk_\s$ is adapted to this decomposition in the sense that 
$$ \fk_\s = \prod_{j \in J} \fk_j \quad \hbox{ for } \quad \fk_j := \fk_\fs \cap \fs_j. $$
Then $\fk' = \fk_\s' \cong \prod_{j \in J} \fk_j'$, where 
$\fk_j'$ is maximal compactly embedded in $\fs_j$. Therefore 
$\fk'$ is finite-dimensional if and only if almost all $\fk_j$ are abelian. 
Since $\fk_j$ is abelian if and only if $\s_j \cong \sL_2(\R)$ (Lemma~\ref{lem:findim-cont}), 
we see that (1) and (3) are equivalent. 

Applying Theorem~\ref{thm:2.5} to the Lie algebra $\fk$, 
we see that its smoothness is equivalent to 
$\fk' \cong \fk/\z(\fk) \cong \fk/\rad(\fk)$ being finite-dimensional, 
which is the equivalence between (2) and (3). 
\end{proof}

\section{Pro-Lie groups as Lie groups} 

In this section we show that for each smooth pro-Lie algebra $\g$, 
the corresponding $1$-connected pro-Lie group $\Gamma(\g)$ carries a regular 
Lie group structure. From that we further derive a description of all 
regular Lie groups $G$ whose Lie algebra is a pro-Lie algebra, as quotients of 
some $\Gamma(\g)$, $\g$ smooth, by some discrete central subgroup $D$. 
In Theorem~\ref{thm:3.8} we further give a very handy characterization of discrete 
central subgroups of $\Gamma(\g)$. 

The Levi decomposition $\g \cong \fr \rtimes \s$ (Theorem~\ref{thm:2.4}) is a 
key tool to obtain the Lie group structure. Let us first recall that 
for the corresponding $1$-connected pro-Lie group $\Gamma(\g)$, we  have 
$$ \Gamma(\g) \cong \Gamma(\fr) \rtimes \Gamma(\s), \quad \hbox{ where } \quad 
\Gamma(\s) \cong \prod_{j \in J} \Gamma(\s_j). $$

Next we show that the coordinates defined by the map $\Phi$ in 
(\ref{eq:3.1}) turn $\Gamma(\fr)$ into a Lie group: 

\begin{proposition} \label{prop:3.1} If $\fr$ is prosolvable, 
$\fn := \oline{[\fr,\fr]}$ and $\fe$ is a closed complement of $\fn$ in $\fr$, 
then $\Gamma(\fr)$ is a Lie group with respect to the manifold structure defined by the 
map $\Phi$ in {\rm(\ref{eq:3.1})}. 
\end{proposition}

\begin{proof} We have to show that multiplication and inversion are smooth maps on 
$\Gamma(\fr)$, which is equivalent to the smoothness of the corresponding maps 
on $\fn \times \fe$. 

Since $\fn$ is pronilpotent, the BCH series defines  a smooth multiplication 
$*$ on $\fn$ satisfying $\exp_{\Gamma(\fr)}(x*y) = \exp_{\Gamma(\fr)}(x)\exp_{\Gamma(\fr)}(y)$ for $x,y \in \fn$ (Remark~\ref{rem:bch-mult}).  
On the other hand, $\fr/\fn$ is abelian. Therefore, 
in the coordinates given by the map $\Phi$, the multiplication takes the form 
$$ (x,y)(x',y') = (x * e^{\ad y}.x' * f(y,y'), y + y'), $$
where $f$ is the map 
$$ f \: \fe \times \fe \to \fn, \quad (y,y') \mapsto 
\exp_{\Gamma(\fn)}^{-1}\Big(\exp_{\Gamma(\fr)}(y)\exp_{\Gamma(\fr)}(y')
\exp_{\Gamma(\fr)}(y+y')^{-1}\Big). $$
Using the description of $\fn$ as a projective limit, we see that the map 
$$ \fe \times \fn \to \fn, \quad (y,x) \mapsto e^{\ad y} x$$ 
is smooth. Therefore it remains to see that $f$ 
is smooth, but this also follows from a straightforward inverse limit argument and 
its validity in all finite-dimensional Lie algebras. 

For the inversion, we obtain from the formula for the product: 
$$ (x,y)^{-1} = (-e^{-\ad y}.(f(y,-y)*x),-y), $$
which implies its smoothness. 
\end{proof}

We are now ready to characterize those pro-Lie algebras 
for which $\Gamma(\g)$ carries a compatible Lie group structure: 

\begin{theorem} \label{thm:3.2} 
For a pro-Lie algebra $\g$, the following are equivalent: 
\begin{description}
\item[\rm(1)] $\Gamma(\g)$ carries a compatible Lie group structure. 
\item[\rm(2)] $\g$ is smooth. 
\end{description}
\end{theorem}

\begin{proof} In view of Theorem~\ref{thm:2.5}, it remains to 
show that (2) implies (1). 

Write $\g \cong \fr \rtimes \s$ with 
$\s \cong \prod_{j \in J} \s_j$ and 
put $$J_0 := \{ j \in J \: \s_j \cong \sL_2(\R)\}, \quad 
\s_0 := \prod_{j \in J_0} \s_j \cong \sL_2(\R)^{J_0}, 
 \quad \mbox{ and } \quad \s_1 := \prod_{j \not\in J_0} \s_j.$$ 
Let $\phi \: \sL_2(\R) \to \tilde\SL_2(\R)$ be a diffeomorphism 
and note that 
$$ \phi^{J_0} \: \s_0  \to \Gamma(\s_0) \cong \tilde\SL_2(\R)^{J_0}$$
defines on $\Gamma(\s_0)$ a Lie group structure (cf.\ Remark~\ref{rem:2.6}). 
Since $\Gamma(\s_1)$ is a finite-dimensional 
Lie group (Theorem~\ref{thm:2.5}), $\Gamma(\s) \cong \Gamma(\s_0) \times \Gamma(\s_1)$ 
is a Lie group. 

In view of $\Gamma(\g) \cong \Gamma(\fr) \rtimes \Gamma(\s)$, it remains to see that 
the action of $\Gamma(\s)$ on $\Gamma(\fr)$ 
is smooth with respect to the Lie group structure on~$\Gamma(\fr)$. 

Clearly, the action of $\Gamma(\s)$ on the Lie algebra $\fn = \oline{[\fr,\fr]}$ 
is smooth, 
because the $\Gamma(\s)$-module 
$\fn$ is a projective limit of finite-dimensional $\Gamma(\s)$-modules on which $\Gamma(\s)$ 
acts smoothly. 
Similarly, it follows that the action of $\Gamma(\s)$ on $\fr/\fn \cong \fe$ is smooth. 

Next we claim that we may choose the closed complement $\fe$ of $\fn$ in 
$\fr$ in an $\s$-invariant fashion. Since $\fr$ is a projective limit 
of finite-dimensional $\s$-modules, its topological dual space is a 
direct limit of finite-dimensional $\s$-modules, hence a semisimple $\s$-module. 
Therefore $\fn^\bot$ has an $\s$-invariant complement $\fm$ in $\fr'$, so 
that we may choose $\fe := \fm^\bot$ (cf.\ \cite{HoMo06}, Theorem~7.16ff). 
Then $\Phi$ is $\Gamma(\s)$-equivariant, which implies the smoothness of 
the $\Gamma(\s)$-action on $\Gamma(\fr)$. 
\end{proof}

The preceding theorem provides one Lie group structure on the topological group 
$\Gamma(\g)$, 
but it is not clear that this is the only one. Therefore it is a 
natural question to ask which additional requirements make this Lie group 
structure unique. One possibility is to require regularity. 
Indeed the regularity of a Lie group $G$, 
which we shall introduce in the next definition,
will secure the existence of an effective exponential function
of $G$ which then allows us to prove the desired uniqueness of
the Lie group structure. 

Let $I$ denote the unit interval $[0,1]$, and abbreviate the Lie algebra
$\L(G)$ of a Lie group $G$ by $\g$.

\begin{definition} A Lie group $G$ is called {\it regular} if for each 
$\xi \in C^\infty(I,\g)$, the initial value problem 
$$ \gamma(0) = \1, \quad \gamma'(t) = \gamma(t)\cdot \xi(t) 
= T_\1(\lambda_{\gamma(t)}) \cdot \xi(t) $$
has a solution $\gamma_\xi \in C^\infty(I,G)$, and the
map 
$$ \eps_G  \: C^\infty(I,\g) \to G, \quad \xi \mapsto \gamma_\xi(1) $$
is smooth (\cite{Mil84}). Then 
$\exp_G(x) := \eps_G(x)$, where $x \in \g$ is identified with a constant function 
$I \to \g$, yields an exponential function of $G$. 
\end{definition}

A crucial feature of regularity is the following (\cite{Mil84}, \cite{GN06}): 
\begin{theorem} \label{thm:reg-homo} If $H$ is a regular Lie group, $G$ is a $1$-connected 
Lie group, and $\phi \: \L(G)\to \L(H)$ is a continuous homomorphism of Lie algebras,
then there exists a unique Lie group homomorphism $f \: G \to H$
with $\L(f) = \phi$. 
\end{theorem}

\begin{corollary} \label{cor:reg-iso} Two 
$1$-connected regular Lie groups with isomorphic Lie algebras 
are isomorphic. 
\end{corollary} 

An important criterion for regularity is provided by the fact that it is an 
extension property. To make this precise, recall that 
an {\it extension of Lie groups} is a short exact sequence 
$$ \1 \to N \sssmapright{\iota} \hat G \sssmapright{q} G \to \1 $$
of Lie group morphisms, for which 
$\hat G$ is a smooth (locally trivial) principal 
$N$-bundle over $G$ with respect to the right action of $N$ given by 
$(\hat g,n) \mapsto \hat gn$, where we identify $N$ with the subgroup 
$\iota(N)$ of $\hat G$.
We call $\hat G$ an {\it extension 
of $G$ by $N$}.

\begin{theorem} \label{thm:reg-ext} {\rm(\cite{KM97}, \cite{GN06})} If 
$\hat G$ is a Lie group extension of $G$ by $N$, 
then $\hat G$ is regular if and only if $G$ and $N$ are regular. 
\end{theorem}

\begin{proposition} \label{prop:reglie} The Lie group structure of $\Gamma(\g)$ from 
{\rm Theorem~\ref{thm:3.2}} turns $\Gamma(\g)$ into a regular Lie group.   
\end{proposition}

\begin{proof} Since $\Gamma(\g)$ is a Lie group extension of $\Gamma(\s)$ by 
$\Gamma(\fr)$ and $\Gamma(\fr)$ a Lie group extension of $(\fe,+)$ by 
the pronilpotent group $(\fn,*)$, it suffices to show that 
$(\fn,*)$, $(\fe,+)$ and $\Gamma(\s)$ are regular. 

In \cite{GN06}, it is shown that for each nilpotent Lie algebra 
$\fm$, the group $(\fm,*)$ is regular, where $*$ denotes the BCH multiplication. 
First, this implies that $(\fe,+)$ is regular. 
Further, $\fn = \prolim \fn_j$ 
for nilpotent Lie algebras $\fn_j$, so that the relation 
$$ C^\infty([0,1],\fn) \cong \prolim C^\infty([0,1],\fn_j) $$
(cf.\ \cite{GN06}) easily implies that $(\fn,*)$ is regular. 

Moreover, all finite-dimensional Lie groups are regular (\cite{KM97}), 
so that a similar argument implies that
$\Gamma(\s) \cong \prod_{j \in J} \Gamma(\s_j)$
is regular. 
\end{proof}

Note that the following theorem deals with Lie groups $G$ whose Lie algebra 
$\L(G)$ is a pro-Lie algebra, but that we do not assume that the underlying 
topological group is pro-Lie.  

\begin{theorem} \label{thm:nec-smooth} A 
pro-Lie algebra $\g$ is the Lie algebra of a 
Lie group with a smooth exponential function if and only if it is smooth. 
\end{theorem}

\begin{proof} Combining Proposition~\ref{prop:reglie} with Theorem~\ref{thm:3.2} 
shows that for each smooth pro-Lie algebra $\g$, the group $\Gamma(\g)$ 
carries a regular Lie group structure with Lie algebra~$\g$, hence in particular 
that the Lie group $\Gamma(\g)$ has a smooth exponential function. 

Assume, conversely, that $G$ is a Lie group with Lie algebra 
$\g$ and a smooth exponential 
function $\exp_G \: \g \to G$. Let $\s = \prod_{j \in J} \s_j$ be a Levi 
complement in $\g$ and $J_0 := \{ j \in J \: \s_j \cong \sL_2(\R)\}.$ 

Since $G$ has a smooth exponential function and each $\s_j$ is 
locally exponential with $\z(\s_j) = \{0\}$, Theorem~IV.4.9 in \cite{Ne06} 
(see \cite{GN06} for a proof) implies the existence of a Lie group morphism 
$\alpha_j \: \Gamma(\s_j) \to G$ for which $\L(\alpha_j)$ is the inclusion map. 
From the regularity of the 
finite-dimensional Lie group $\Gamma(\s_j)$, we further obtain with Theorem~\ref{thm:reg-homo} 
morphisms $\beta_j \: G \to \Gamma(\s_j)$, for which $\L(\beta_j)$ is the 
projection $\g \to \s_j$. We then have 
$\beta_j \circ \alpha_j = \id_{\Gamma(\s_j)}$ for each $j \in J$. 

Let $U \subeq G$ be a contractible $\1$-neighborhood and 
$\fk_\s \subeq \s$ a maximal compactly embedded subalgebra. 
Then $\fk_\s$ is adapted to the decomposition of $\s$ in the sense that 
$\fk_\s = \prod_{j \in J} \fk_j$ for $\fk_j := \fk_\fs \cap \fs_j.$

Since the exponential function of $G$ is continuous, $V := \exp_G^{-1}(U)$ 
is a $0$-neighborhood in $\g$, and we conclude that 
$\fk_\s \cap V$ contains a subalgebra of the form 
$\prod_{j \in J_1} \fk_j,$
where $J_1 \subeq J$ is a cofinite subset. 

For each $j \in J$, the group 
$\exp_{\Gamma(\s_j)}(\fk_j)$ is a maximal compactly embedded subgroup 
$K_j$ of $\Gamma(\s_j)$ and the inclusion $K_j \into \Gamma(\s_j)$ is a 
homotopy equivalence 
by \cite{Ho65}, p.~180, Theorem~3.1. For $j \in J_1$, we have 
$$\alpha_j(\exp_{\Gamma(\s_j)}\fk_j) 
= \exp_G(\fk_j) \subeq U, $$
so that the contractibility of $U$ in $G$ implies that 
the map $\alpha_j \: \Gamma(\s_j) \to G$ is homotopic to a constant map, 
hence that $\id_{\Gamma(\s_j)} = \beta_j \circ \alpha_j$ is also 
homotopic to a constant map, that is,  $\Gamma(\s_j)$ is contractible. 
In view of Lemma~\ref{lem:findim-cont}, this implies 
$\s_j \cong \sL_2(\R)$, hence that $J_1 \subeq J_0$, so that 
$J_0$ is cofinite, that is,  $\g$ is smooth (Theorem~\ref{thm:2.5}). 
\end{proof}

\begin{corollary} \label{cor:reg-cover} If $G$ is a regular Lie group for which 
$\g = \L(G)$ is a pro-Lie algebra, then $\g$ is smooth and 
$G$ is isomorphic to a quotient 
of the regular Lie group $\Gamma(\g)$ by a discrete central subgroup. 
In particular, $G$ is a pro-Lie group. 
\end{corollary}

\begin{proof} If $G$ is any connected 
regular Lie group whose Lie algebra $\L(G)$ is a pro-Lie algebra, 
then its universal covering group $\tilde G$ is a regular $1$-connected Lie group 
with Lie algebra $\g$, and Theorem~\ref{thm:nec-smooth} implies that 
$\g$ is smooth. Hence Proposition~\ref{prop:reglie} 
and Corollary~\ref{cor:reg-iso} imply that $\tilde G \cong \Gamma(\g)$. 
Now the assertion follows from the fact that the universal covering map 
$q_G \: \tilde G \to G$ has discrete central kernel 
(see \cite{HoMo06}, Lemma~3.32(ii)). 
\end{proof}

In view of the preceding corollary, it is 
of crucial importance to understand the discrete central 
subgroups of the groups $\Gamma(\g)$, provided 
this group carries a Lie group structure. 
  
\begin{lemma} \label{lem:3.5} 
Any discrete central subgroup $\Gamma$ of a pro-Lie group $G$ is finitely generated. 
\end{lemma}

\begin{proof} Let $U \subeq G$ be an open identity neighborhood with $U \cap \Gamma = \{\1\}$. 
After shrinking $U$, if necessary, we may further assume that there exists a 
closed normal subgroup $N \trile G$ with $NU = UN = U$ such that 
$G/N$ is a finite-dimensional Lie group (Theorem~\ref{thm:2.0b}). 
Then $U/N$ is an open identity neighborhood of $G/N$ intersecting $\Gamma N/N$ 
trivially. Hence $\Gamma \cong \Gamma N/N$ is a discrete central subgroup of 
$G/N$, and therefore finitely generated.\begin{footnote}{To see that any discrete 
central subgroup $\Gamma$ of a connected finite-dimensional Lie group $G$ 
is finitely generated, we first recall from \cite{Ho65} that $Z(G)$ is contained 
in a connected abelian Lie subgroup $A$ of $G$. Since $\tilde A \cong \R^n$ 
for some $n$, it suffices to observe that discrete subgroups of $\R^n$ 
are finitely generated because they are isomorphic to $\Z^m$ for some $m \leq n$.}
\end{footnote}
\end{proof}

For an abelian topological group $A$, we write $\comp(A)$ for the 
subgroup, generated by all compact subgroups of $A$. 

\begin{lemma} \label{lem:3.6} 
If $\g$ is smooth and $Z := Z(\Gamma(\g))$, then 
$\comp(Z) = \tor(Z)$ 
is a finite group. 
\end{lemma}

\begin{proof} We have seen in (\ref{eq:3.1}) that $\Gamma(\fr)$ is an extension 
of the additive group $(\fe, +)$ by the $1$-connected pronilpotent group 
$(\fn,*)$. Since both these groups are compact free, $\Gamma(\fr)$ is compact free. 

As in the proof of Theorem~\ref{thm:3.2}, we write $\s = \s_0 \times \s_1$, where 
$\s_0 \cong \sL_2(\R)^{J_0}$ and $\s_1$ is finite-dimensional semisimple. Then 
\begin{eqnarray}
  \label{eq:3.2} 
 Z(\Gamma(\s)) = Z(\Gamma(\s_0)) \times 
Z(\Gamma(\s_1)) \cong \Z^{J_0} \times Z(\Gamma(\s_1)),
\end{eqnarray}
where $Z(\Gamma(\s_1))$ is a finitely generated abelian group. 

If $C \subeq Z$ is a compact subgroup, then $C \cap \Gamma(\fr) = \{\1\}$ 
($\Gamma(\fr)$ is compact free) implies that $C$ injects into $Z(\Gamma(\s))$, and 
since $Z(\Gamma(\s_0))$ is compact free, 
$C$ injects into $Z(\Gamma(\s_1))$. As $Z(\Gamma(\s_1))$ is finitely generated, 
$\comp(Z(\Gamma(\s_1))) = \tor(Z(\Gamma(\s_1)))$ is a finite group. We conclude that 
$C$ consists of torsion elements, which already implies that 
$\comp(Z) = \tor(Z)$. As $\tor(Z)$ intersects $\Gamma(\fr)$ trivially, 
it also injects into $Z(\Gamma(\s))$, hence into $Z(\Gamma(\s_1))$, which implies 
its finiteness. 
\end{proof}

\begin{lemma} \label{lem:3.7} 
For an arbitrary set $J$, a finitely generated subgroup of $\Z^J$ is discrete. 
\end{lemma}

\begin{proof} Let $\Gamma \subeq \Z^J$ be a finitely generated subgroup. 
Then $\Gamma$ is torsion free, hence isomorphic to $\Z^d$ for some $d \in \N_0$. 
Let $\chi_j \: \Z^J \to \Z$, $j \in J$, denote the coordinate projections. 
Then the restrictions $\oline\chi_j \: \Gamma \to \Z$ separate the points, 
hence generate a subgroup of full rank in $\Hom(\Gamma,\Z) \cong \Z^d$. 
We conclude that there exists a finite subset $F \subeq J$ such that 
the kernel of the projection 
$$\chi_F := (\chi_j)_{j \in F} \: \Z^J \to \Z^F $$ 
intersects $\Gamma$ trivially. Since $\ker(\chi_F)$ is an open subgroup, 
$\Gamma$ is discrete. 
\end{proof}

The following theorem provides a very nice characterization of the 
discrete central subgroups of $\Gamma(\g)$ for a smooth pro-Lie algebra 
$\g$ and hence a description of all regular Lie groups $G$ with pro-Lie algebras 
as Lie algebras~$\L(G)$. 

\begin{theorem} \label{thm:3.8} Assume that $\g$ is smooth. 
Then a subgroup $\Gamma$ of $Z := Z(\Gamma(\g))$ is discrete if 
and only if it is finitely generated and $\Gamma \cap Z_0$ is discrete.   
\end{theorem} 

\begin{proof} In view of Lemma~\ref{lem:3.5}, each discrete subgroup of 
$Z$ is finitely generated. 

Conversely, assume that $\Gamma \subeq Z$ is finitely generated and that 
$\Gamma \cap Z_0$ is discrete. Then Theorem~5.32(iv) in \cite{HoMo06} 
implies that the pro-Lie group $\oline\Gamma$ (cf.\ Theorem~\ref{thm:2.0}) 
is a direct product $\oline\Gamma \cong \R^m\times \comp(\oline\Gamma)\times\Z^n$ for 
some $m,n\in\N_0$. Since $\comp(Z)$ is finite (Lemma~\ref{lem:3.6}), 
$\pi_0(\oline\Gamma)$ is discrete and $\oline\Gamma_0 \cong \R^m$. 
As $\oline\Gamma_0$ is contained in $Z_0 \cong \z(\g)$, it is contained in 
the closure of $\Gamma \cap Z_0$, which is discrete by assumption, hence closed. Since 
$\Gamma \cap Z_0$ is countable, we get $m = 0$, so that $\Gamma$ is discrete. 
\end{proof}

\section{Locally exponential Lie algebras} 

In this section, we recall some basic definitions and properties 
concerning locally exponential Lie algebras. In particular, we 
introduce the Maurer--Cartan form and derive a spectral 
condition from the invertibility properties of the Maurer--Cartan form. 
This section prepares the following one, where we characterize the locally exponential 
pro-Lie algebras, respectively, pro-Lie groups.

\begin{definition} A Lie group $G$ is called {\it locally exponential} 
if it has a smooth exponential function $\exp_G \: \L(G) \to G$ 
mapping some open $0$-neighborhood in $\L(G)$ diffeomorphically 
onto an open $\1$-neighborhood in~$G$.
\end{definition}

\begin{remark} \label{rem:4.3} (a) The Lie algebra $\L(G)$ of a locally exponential 
Lie group $G$ is locally exponential (\cite{GN06}; \cite{Ne06}, Lemma~IV.2.2). 

(b) All Banach--Lie algebras and therefore all finite-dimensional 
Lie algebras are locally exponential because the BCH series defines a 
smooth local group structure on some $0$-neighborhood in $\g$, 
satisfying all requirements of Definition~\ref{def:2.2}. 
\end{remark}

\begin{lemma} {\rm(\cite{GN06})} If 
$\g$ is locally exponential, then all operators 
$\ad x$ generate a smooth $\R$-action 
$(t,y) \mapsto e^{t\ad x}y$ on $\g$ by automorphisms of topological 
Lie algebras. 
\end{lemma}

\begin{definition}\label{def:4.5} 
Now let $\g$ be locally exponential and $U$ as in Definition~\ref{def:2.2}. 
Then, for each $x \in U$, the left multiplication $\lambda_{-x}^* \: y \mapsto (-x) * y$ 
is defined in a neighborhood of $x$ with $\lambda_{-x}(x) = (-x) * x  =0$, and 
$$ \kappa_U(x) := T_x(\lambda_{-x}) \: \g \to \g $$
defines a $\g$-valued $1$-form on $U$, called the {\it Maurer--Cartan form}.
In \cite{GN06} it is shown that the Maurer--Cartan form can be expressed 
by an operator-valued integral: 
$$ (\kappa_U)_x = \kappa_\g(x) := \int_0^1 e^{-t \ad x}\, dt, $$
interpreted in the pointwise sense. If $\g$ is complete, 
the integral $\kappa_\g(x)$ 
is defined for each $x\in \g$, but for $x \in U$, the interpretation 
in terms of the Maurer--Cartan form implies that the linear operator 
$\kappa_\g(x)$ is 
invertible. 

We call a point $x \in \g$ {\it $\exp$-regular} if the operator $\kappa_\g(x)$ 
is invertible. This terminology is justified by the fact that if 
$\exp_G \: \L(G) \to G$ is the exponential function of a Lie group, 
then $T_x(\exp_G)$ is invertible if and only if $x$ is $\exp$-regular. 
\end{definition}

\begin{remark} \label{rem:4.6} If $\g$ is a pro-Lie algebra, we write 
$\g = \prolim \g_j$ with $\g_j \cong \g/\fn_j$ for closed ideals 
$\fn_j \trile \g$ of finite codimension (\cite{HoMo06}, Definition~3.6 and Proposition~3.9). 
Then $\kappa_\g(x)$ preserves 
each ideal $\fn_j$ and induces an operator on the quotient 
$\g/\fn_j$. Since all these quotients are finite-dimensional, 
$\kappa_\g(x)$ is invertible if and only if all operators induced 
on the quotients $\g/\fn_j$ are invertible, which is equivalent to 
\begin{eqnarray}
  \label{eq:reg-cond}
\Spec(\ad x) \cap 2\pi i \Z = \{0\}. 
\end{eqnarray}
\end{remark}

\begin{remark} If $G$ is  a locally exponential Lie group and 
$\L(G)$ its Lie algebra, then the exponential function induces a 
one-to-one map $\L(G) \to \fL(G), x \mapsto \gamma_x$ which is also 
a homeomorphism (\cite{GN06}). If, in addition, $G$ is a pro-Lie group, 
we conclude that $\L(G) \cong \fL(G)$ are topological Lie algebras, 
and hence that $\g := \fL(G)$ is a locally exponential pro-Lie algebra. 
This further implies that $\tilde G \cong \Gamma(\g)$ 
(\cite{HoMo06}, Theorem~8.21). Hence 
$G \cong \Gamma(\g)/D$ for some discrete central subgroup $D$ of $\Gamma(\g)$. 

If, conversely, $\g$ is a locally exponential pro-Lie algebra, 
then $\Gamma(\g)$ is a locally exponential Lie group and for each discrete 
central subgroup $D \subeq \Gamma(\g)$, the quotient 
$\Gamma(\g)/D$ is a locally exponential Lie group with Lie algebra~$\g$. 
\end{remark}

\section{Locally exponential pro-Lie algebras} 

The main result of this section is a characterization of locally 
exponential pro-Lie algebras as those for which the $\exp$-regular points 
form a $0$-neighborhood. Moreover, we shall see that this condition implies 
that $\g$ is smooth and that the corresponding simply connected regular 
Lie group $\Gamma(\g)$ is locally exponential. 

\subsection{Exponential pro-Lie algebras} 

We start with a discussion of exponential Lie algebras. 

\begin{theorem} \label{thm:4.1} For a pro-Lie algebra $\g$ the following are 
equivalent: 
\begin{description}
\item[\rm(1)] $\g$ satisfies the spectral condition 
  \begin{description}
    \item[\rm(SC)] \qquad $(\forall x \in \g) \quad \Spec(\ad x) \cap i \R = \{0\}.$ 
  \end{description}
\item[\rm(2)] $\Gamma(\g)$ is an exponential Lie group.  
\item[\rm(3)] $\g$ is exponential.  
\end{description}
Any such Lie algebra is prosolvable. 
\end{theorem}

\begin{proof} (1) $\Rarrow$ (2): Assume the spectral condition (SC). We observe that 
it implies that for each closed ideal $\fn \trile \g$ and the quotient map 
$q \: \g \to \fq := \g/\fn$ we have 
$$\Spec(\ad_\fq(q(x))) \cap i \R \subeq \Spec(\ad_\g(x))\cap i \R= \{0\}. $$
Hence $\fq$ is exponential and therefore solvable by the Dixmier--Saito Theorem 
(\cite{Dix57}, \cite{Sai57}). 

We thus have $\g \cong \prolim \g_j$, where each $\g_j$ is exponential, so that 
$\g$ is in particular prosolvable. 

If $\phi_{ij} \: \g_j \to \g_i$ is a homomorphism of exponential Lie algebras, 
then the corresponding homomorphism $\tilde\phi_{ij} \: (\g_j,*) \to (\g_i,*)$ 
of simply connected Lie groups satisfies 
$$ \tilde\phi_{ij} \circ \exp_{(\g_j,*)} = \exp_{(\g_i,*)} \circ \phi_{ij}, $$
which shows that  $\tilde\phi_{ij} = \phi_{ij}$ also respects the 
$*$-product, so that we obtain a $*$-product on any projective 
limit $\g = \prolim \g_j$ of exponential Lie algebras, showing that 
$\g$ is exponential. This implies in particular that $\Gamma(\g) \cong (\g,*)$ is an 
exponential Lie group. 

(2) $\Rarrow$ (3) follows from $\g \cong \L(\Gamma(\g)) \cong \fL(\Gamma(\g))$ 
if $\Gamma(\g)$ is an exponential Lie group. 

(3) $\Rarrow$ (1): To see that (SC) is satisfied if 
$\g$ is exponential, we simply observe that if 
$\g$ is exponential, the operator $\kappa_\g(x) = \int_0^1 e^{-t\ad x}\, dt$ 
is invertible for each $x \in \g$, that is,  $\Spec(\ad x) \cap 2\pi i \Z \subeq \{0\}$, 
and this implies (SC) (Remark~\ref{rem:4.6}).  
\end{proof}

For finite-dimensional solvable Lie algebras, it is quite convenient 
to have Saito's testing device (\cite{Sai57}, see also \cite{Bou89}, Ch.~III, Ex.~9.17), characterizing the exponential 
Lie algebras as the solvable Lie algebras 
not containing a subalgebra isomorphic to 
$\mot_2$, the Lie algebras of the motion group of the euclidean plane, or 
its four-dimensional central extension $\osc$, the oscillator algebra. 
These Lie algebras can be described in terms of commutator relations 
as follows. The $3$-dimensional Lie algebra $\mot_2$ has a basis 
$U,P,Q$ with 
$$ [U,P] = Q, \quad [U,Q] = -P \quad \mbox{ and } \quad [P,Q] = 0, $$
whereas $\osc$ has a basis $U,P,Q,Z$, where $Z$ is central with 
$$ [U,P] = Q, \quad [U,Q] = -P \quad \mbox{ and } \quad [P,Q] = Z. $$
One implication of Saito's result is trivial, because in both Lie algebras 
we have $i \in \Spec(U)$, so that the occurence of any such subalgebra in a 
Lie algebra $\g$ implies that $\g$ is not exponential. 

Conversely, any finite-dimensional non-exponential Lie algebra $\g$ contains 
a triple $(U,P,Q)$ satisfying 
\begin{eqnarray} \label{com-rel} 
[U,P] = Q \quad \mbox{ and } \quad [U,Q] = -P.   
\end{eqnarray}
One finds more such pairs as follows: Put 
$P_1 := P$, $Q_1 := Q$ and $Z_1 := [P,Q]$ and, recursively, 
$P_{i+1} := [Z_i, P_i]$, $Q_{i+1} := [Z_i, Q_i]$, $Z_{i+1} := [P_{i+1}, Q_{i+1}]$. 
If $\g$ is solvable, then it is easy to see that 
for some $n$ the elements $U, P_n, Q_n, Z_n$ span a subalgebra either isomorphic 
to $\mot_2$ or $\osc$ (\cite{Sai57}). Here the main point is that 
$Z_i \in D^i(\g)$ vanishes if $i$ is large enough. 

One might expect that similar testing devices exist for prosolvable 
Lie algebras. However, the following example shows 
that the situation becomes more complicated. 

\begin{example} Choose a basis 
$$ P := \frac{1}{2}\begin{pmatrix} 0 & 1 \cr 1 & 0 \end{pmatrix}, \quad 
Q := \frac{1}{2}\begin{pmatrix}1 & 0 \cr 0 & -1 \end{pmatrix}, \quad 
U := \frac{1}{2}\begin{pmatrix}0 & 1 \cr -1 & 0 \end{pmatrix} $$
in $\sL_2(\R)$, which satisfies the commutator relations 
$$ [U,P] = Q, \quad [U,Q] = -P \quad \mbox{ and } \quad [P,Q] = U. $$

Let $A := X\R[[X]]$ be the pronilpotent algebra of formal power series 
in one variable vanishing in $0$. We observe that $A$ is formally real in the sense 
that for $(f,g) \not=(0,0)$ in $A^2$ we have $f^2 + g^2 \not=0$. 
Then $A \otimes \sL_2(\R) \cong \sL_2(A)$ is a pronilpotent 
Lie algebra with respect to the bracket defined by 
$$ [a \otimes x, a' \otimes x'] = aa' \otimes [x,x']. $$
Define 
$$\g := (A \otimes \sL_2(\R)) \rtimes (\R \otimes U) 
\subeq \R[[X]] \otimes \sL_2(\R) \cong 
\sL_2(\R[[X]])$$ 
and note that $\g \cong \sL_2(A) \rtimes \R U$ 
is a prosolvable Lie algebra with a pronilpotent hyperplane ideal 
and $i \in \Spec(\ad U)$, so that $\g$ is not exponential. 
For any nonzero pair $P,Q \in \g$ satisfying (\ref{com-rel}), we then have 
$$ P = a \otimes P + b \otimes Q \quad \mbox{ and }\quad  Q = -b \otimes P + a \otimes Q, $$
and therefore 
$$ Z := [P,Q] = (a^2 + b^2) \otimes U \not=0. $$
Moreover, $[Z,P]$ and $[Z,Q]$ are nonzero, so that the recursive 
construction from above produces infinitely many nonzero elements and 
Saito's method to find subalgebras isomorphic to $\mot_2$ or $\osc$ breaks down. 
We also note that for $0\not=c \in X\R[[X]]$ the element  
$U' := (1+ c) \otimes U \in \g$ satisfies $i \in \Spec(U')$ but the operator 
$(\ad U')^2 + \1$ is injective. 
\end{example}


\subsection{Locally exponential pro-Lie algebras} 

As we shall see below, it requires some work
 to characterize the locally exponential pro-Lie algebras, 
but it is easy to find a strong necessary condition which already 
provides the key hint on how to approach the problem. 

\begin{proposition} \label{prop:4.2} 
If the pro-Lie algebra $\g$ is 
locally exponential, then it contains a closed exponential ideal of finite codimension. 
In particular, $\g/\rad(\g)$ is finite-dimensional. 
\end{proposition}

\begin{proof} If $\g$ is locally exponential, then Definition~\ref{def:4.5} 
implies the existence of a $0$-neighborhood $U \subeq \g$ 
such that $\kappa_\g(x)$ is invertible for each $x \in U$. 
Since $\g$ is a projective limit of finite-dimensional Lie algebras, 
$U$ contains a closed ideal $\fn$ of finite codimension. 
From Remark~\ref{rem:4.6} we further see that this implies that 
$\Spec(\ad x) \cap 2 \pi i \Z = \{0\}$ for each $x \in U$. 

For each $x \in \fn$, we derive from 
$$ \Spec(\ad_\g(x)) = \Spec(\ad_\fn(x)) \cup \{0\} $$ 
and $\R x \subeq U$ 
that $\Spec(\ad_\fn tx) \cap 2\pi i \Z = \{0\}$ for each $t \in \R$, 
so that $\Spec(\ad_\fn(x)) \cap i \R = \{0\}$, and Theorem~\ref{thm:4.1} 
implies that $\fn$ is exponential and hence prosolvable. 
\end{proof}

\begin{remark} If $\g$ has a closed exponential ideal 
of finite codimension, 
then $\g/\rad(\g)$ is finite-dimensional, which implies in particular that $\g$ is smooth 
and hence that $\Gamma(\g)$ carries 
a regular Lie group structure (Proposition~\ref{prop:reglie}). 
\end{remark}

Next we describe an example of a pro-Lie algebra with an abelian 
hyperplane ideal which is not locally exponential, so that the condition 
in Proposition~\ref{prop:4.2} is not sufficient for local exponentiality. 

\begin{example} \label{ex:4.3} 
Let $\alpha \: \R \to \GL(E)$ be a smooth 
representation of $\R$ on a  
complete locally convex space $E$ with the 
infinitesimal generator $D = \alpha'(0)$. 
Then the semi-direct product group 
$$ G := E \rtimes_\alpha \R, \quad (v,t) (v',t') = (v + \alpha(t)v', t + t') $$
is a Lie group with Lie algebra $\g = E \rtimes_D \R$ and exponential function 
$$ \exp_G(v,t) 
= \big( \beta(t)v,  t\big), \qquad 
\beta(t) =  \int_0^1 \alpha(st)\, ds = 
\left\{ \begin{array}{cl} 
\id_E & \mbox{ for $t = 0$} \cr 
\frac{1}{t}\int_0^t \alpha(s)\, ds & \mbox{ for $t \not=0$.}
\end{array} \right.$$ 
From this formula it follows that $(w,t) \in \im(\exp_G)$ is equivalent to 
$w \in \im(\beta(t))$. We conclude that 
$\exp_G$ is injective on some $0$-neighborhood 
if and only if $\beta(t)$ is injective for $t$ close to $0$, 
and it is surjective onto some $\1$-neighborhood in $G$ if 
and only if $\beta(t)$ is surjective for $t$ close to $0$. 

Note that the eigenvector equation $Dv = \lambda v$ for $t\lambda \not=0$ implies that 
$$ \beta(t)v = \int_0^1 e^{st\lambda}v\, ds  = \frac{e^{t\lambda} - 1}{t\lambda}v, $$
so that $\beta(t)v = 0$ is equivalent to $t\lambda \in 2 \pi i \Z \setminus \{0\}$. 

(a) For the weakly complete space $E = \C^\N$ 
and the diagonal operator $D$ given by 
$D((z_n)_{n \in \N}) = (2\pi in z_n)_{n \in \N}$, we see that 
$\beta(\frac{1}{n}) e_n = 0$ holds for $e_n = (\delta_{mn})_{m \in \N}$, 
and $e_n \not\in \im\big(\beta(\frac{1}{n})\big)$. 
We conclude that $(e_n, \frac{1}{n})$ is not contained in the image of $\exp_G$, and since 
$(e_n, \frac{1}{n})\to (0,0)$, the identity of $G$, 
$\im(\exp_G)$ does not contain any identity neighborhood of $G$. 
Hence the exponential function of the Lie group 
$G = E \rtimes_\alpha \R$ is neither locally injective nor locally surjective in $0$. 

(b) For the Fr\'echet space $E = \R^\N$ 
and the diagonal operator $D$ given by 
$D((z_n)_{n \in \N}) = (n z_n)_{n \in \N}$, 
it is easy to see that all operators $\beta(t)$ are invertible 
and that $\R \times E \to E \times E, (t,v) \mapsto (\beta(t)v, \beta(t)^{-1}v)$ 
is a smooth map. This implies that $\exp_G \: \g \to G$ is a diffeomorphism. 
\end{example}

\begin{definition} \label{def:root} Let $\g$ be a prosolvable Lie algebra. 
A {\it root of $\g$} is a continuous linear functional 
$\alpha \: \g_\C \to \C$ with the property that there exist two closed 
ideals $\fn_1 \subeq \fn_2$ of $\g_\C$ with $\dim_\C (\fn_2/\fn_1) = 1$ such that 
$$ (\forall x \in \g)\quad (\ad x - \alpha(x)\1)(\fn_2) \subeq \fn_1. $$
We write $\Delta(\g) \subeq \Hom_\C(\g_\C,\C)$ for the set of roots of $\g$. 
\end{definition}

\begin{lemma} The roots of a prosolvable Lie algebra $\g$ vanish on the commutator 
subalgebra, hence can be interpreted as homomorphisms of Lie algebras. 
For each element $x \in \g$, we have 
$$ \Spec(\ad x) = \Delta(\g)(x). $$
\end{lemma}

\begin{proof} Let $\alpha \in \Delta(\g)$, 
choose $\fn_1$ and $\fn_2$ as in the definition, then 
we have a representation 
 $\delta$ of $\g$ into $\gl(\fn_2/\fn_1)$, given by
$\delta(x)(v+\fn_1)=[x,v]+\fn_1$. Since $\fn_2/\fn_1$ is
one-dimensional, the Lie algebra $\gl(\fn_2/\fn_1)$ is abelian, and thus
$\delta(\g')=\{0\}$. But by Definition~\ref{def:root}, we have
$\delta(x)(v+\fn_1)=\alpha(x)\cdot (v+\fn_1)$, and thus $\alpha(\g')=\{0\}$.

For each $x \in \g$ and each root $\alpha$, we clearly have $\alpha(x) \in \Spec(\ad x)$. 
Conversely, each spectral value $\lambda$ of $\ad x$ is contained in 
$\Spec(\ad_\fq q(x))$ for some quotient map $q \: \g \to \fq$ onto some finite-dimensional 
Lie algebra $\fq$. 
Applying Lie's Theorem to the finite-dimensional solvable Lie algebra $\fq_\C$, 
we see that 
there exist ideals $\fn_1^\fq \subeq \fn_2^\fq$ and a root 
$\alpha^\fq$ of $\fq$ with $\alpha^\fq(q(x)) = \lambda$. 
Then $\alpha := \alpha^\fq\circ q$ is a root of $\g$ with 
$\alpha(x) = \lambda$.  
\end{proof}

\begin{lemma} \label{lem:4.6} Let $\g$ be a finite-dimensional solvable Lie algebra 
and $\Delta(\g) \subeq \Hom(\g_\C,\C)$ its set of roots. 
If $g \in \Gamma(\g)$ satisfies 
$$ \Gamma(\alpha)(g)\not\in 2\pi i \Z \setminus \{0\} 
\quad \mbox{ for each} \quad \alpha \in \Delta(\g), $$ 
then there exists a unique $x \in \g$ with $g = \exp_G(x)$. 
\end{lemma}

\begin{proof} This lemma follows from Theorem 2 in \cite{Dix57}. We reproduce the 
argument for the sake of completeness.  
Put $G := \Gamma(\g)$. 
First we show that $x$ is unique. If $x,y \in \g$ satisfy 
$\exp_G(x)= \exp_G(y) = g$, then the assumption on $g$ implies that 
$x$ is $\exp$-regular, so that $[x,y] = 0$ and $\exp_G(x-y) = \1$ 
(\cite{HHL89}, Lemma~V.6.7; see also \cite{GN06} for the infinite-dimensional case). 
Since $G$ is simply connected and solvable, we get $x - y = 0$, 
because all compact subgroups of $G$ are trivial (\cite{Ho65}). 

Let $D^k(\g) \trile \g$ be the last nonzero 
term of the derived series of $\g$ and $q \: \g \to \fq := \g/D^k(\g)$ the quotient map. 
Then $D^k(\g) \trile \g$ is an abelian ideal of $\g$  
and $\fq$ is a solvable Lie algebra of length $k$, 
whereas the length of $\g$ is~$k+1$. 

If $\g$ is abelian, then $G \cong (\g,+)$ implies that the 
exponential function of $G$ is surjective. We now argue by induction on 
the solvable length of $\g$. The abelian case corresponds to solvable length~$\leq 1$. 
As $\g$ is finite-dimensional, there exists a $1$-connected Lie group 
$Q$ with $\L(Q) = \fq$ and a quotient morphism 
$q_G \: G \to Q$ of Lie groups. 
As $Q$ is $1$-connected, the subgroup $D^k(G) := \ker q_G$ is connected. 

Since all characters of $G$ vanish on $D^1(G)$, they vanish in particular 
on $D^k(G)$, so that each character $\chi$ factors through a homomorphism 
$\chi_Q \: Q \to \C$ with $\chi_Q \circ q_G = \chi$. 
As $\chi_Q(q_G(g))$ never is a nonzero integral multiple of 
$2\pi i$, the same holds for the values of the characters of $Q$ on $q_G(g)$.  
We now apply our induction hypothesis to $Q$ to derive that 
$q_G(g) \in \exp_Q(\fq)$. 
This means that there exists some $x \in \g$ with 
$g \in \exp_G(x) D^k(G)$, so that it remains to show that 
the subset $\exp_G(x) D^k(G)$ of $G$ is contained in the image of 
the exponential function. 

The group $D^k(G)$ is abelian and connected, so that 
its exponential function is surjective. 
We may therefore assume that $g \not\in D^k(G)$, which implies that 
$x \not\in D^k(\g)$. Then $\fb := D^k(\g) + \R x$ is a subalgebra of $\g$, 
isomorphic to $D^k(\g) \rtimes_D \R$, where $D := \ad x\res_{D^k(\g)}$. 
We conclude that with $\alpha(t) := e^{tD}$ we obtain a simply connected group 
$B := D^k(\g) \rtimes_\alpha \R$
with Lie algebra $\fb$. The exponential function of $B$ is given by 
$$ \exp_B(v,t) = (\beta(t)v, t), 
\quad \mbox{ where } \quad \beta(t) = \int_0^1 e^{-st D}\,ds 
= \frac{\1 - e^{-t\ad x}}{t \ad x}.$$
From the spectral condition on $g$ we conclude that 
the operator $\beta(1)$ is invertible, 
which immediately implies that $D^k(\g) \times \{1\} \subeq \im(\exp_B)$. 

Now let $j_B \: B \to G$ be the unique Lie group homomorphism for which 
$\L(j_B) \: \fb \to \g$ is the natural inclusion map. 
Then 
$$g\in D^k(G) \exp_G(x) = j_B(D^k(\g) \times \{1\}) \subeq j_B(\im(\exp_B)) = \exp_G(\fb) 
$$ 
implies that $g \in \exp_G(\g)$. 
\end{proof}

\begin{proposition} \label{prop:4.7} Let $\g$ be a prosolvable pro-Lie algebra 
and $\Delta(\g)$ its set of roots. 
If $g \in \Gamma(\g)$ satisfies 
\begin{eqnarray}
  \label{eq:sc}
\Gamma(\alpha)(g) \not\in 2\pi i \Z \setminus \{0\} \quad \mbox{ for each} \quad \alpha 
\in \Delta(\g), 
\end{eqnarray}
then there exists a unique $x \in \g$ with $g = \exp(x)$. 
\end{proposition}

\begin{proof} We write $\g = \prolim \g_j$ for a family of finite-dimensional 
solvable Lie algebras $\g_j$. We may assume that the corresponding maps 
$q_j \: \g \to \g_j$ are quotient maps (\cite{HoMo06}, Theorem~A2.12). 

On the group level, by 
\cite{HoMo06}, Theorem~6.8, 
we obtain quotient maps 
$\Gamma(q_j) \: \Gamma(\g) \to \Gamma(\g_j)$
of the corresponding $1$-connected Lie groups. 
Viewing $\g_j$ as a $\g$-module, we can think  of the 
roots of $\g_j$ as obtained by factorization of certain roots of $\g$. 
Hence $g_j := \Gamma(q_j)(g)$ satisfies for each $j$ condition (\ref{eq:sc}). 
In view of Lemma~\ref{lem:4.6}, there exists a unique $x_j \in \g_j$ with 
$\exp_{\Gamma(\g_j)}(x_j) = g_j$, and the uniqueness implies that the family 
$(x_j)_{j \in J} \in \prod_{j \in J} \g_j$ 
defines an element $x$ of $\g = \prolim \g_j$ with 
$\exp_{\Gamma(\g)} x = g$. The uniqueness assertion of Lemma~\ref{lem:4.6} also implies 
the uniqueness of $x$. 
\end{proof}

\begin{proposition} \label{prop:cont-smooth} Let 
$f \: G \to H$ be a continuous homomorphism of regular Lie groups 
whose Lie algebras $\g$, respectively, $\h$ are pro-Lie algebras. 
Then $f$ is smooth. 
\end{proposition}

\begin{proof} Since smoothness is a local property, we may without loss of generality  
assume that $G$ and $H$ are $1$-connected; otherwise we replace them 
by the simply connected covering of their identity component 
and $f$ by the induced homomorphism of these $1$-connected groups. 

Now Corollary~\ref{cor:reg-cover} implies that 
$G \cong \Gamma(\g)$ and $H \cong \Gamma(\h)$ are regular Lie groups. 
The homomorphism $f$ is uniquely determined 
by the relation 
$$ f \circ \exp_{\Gamma(\g)} = \exp_{\Gamma(\h)} \circ \fL(f). $$
Since the group $\Gamma(\h)$ is a regular Lie group, 
the continuous 
homomorphism of Lie algebras $\fL(f) \: \g \to \h$ integrates to a unique 
smooth morphism of Lie groups 
$h \: \Gamma(\g) \to \Gamma(\h)$ (\cite{Mil84}) with $\L(h) = \fL(f)$, also satisfying 
$$ h \circ \exp_{\Gamma(\g)} = \exp_{\Gamma(\h)} \circ \fL(f). $$
This implies that $f = h$ and hence that $f$ is smooth. 
\end{proof}

\begin{theorem} \label{thm:locexp} For a pro-Lie algebra $\g$, the following are equivalent: 
  \begin{description}
  \item[\rm(1)] $\g$ is locally exponential. 
  \item[\rm(2)] There exists a $0$-neighborhood 
$U \subeq \g$ consisting of $\exp$-regular points. 
  \item[\rm(3)] $\Gamma(\g)$ is a locally exponential Lie group. 
  \end{description}
\end{theorem}

\begin{proof} (1) $\Rarrow$ (2) is a direct consequence of the discussion in 
Remark~\ref{rem:4.6}. 

(2) $\Rarrow$ (3): Let $\fn \trile \g$ be a closed ideal of finite-codimension 
for which $U$ contains an $\fn$-saturated  $0$-neighborhood and recall from 
Proposition~\ref{prop:4.2} that $\fn$ is exponential and prosolvable.  
We consider the regular Lie group $G := \Gamma(\g)$ with Lie algebra $\g$. 

Let $\fe \subeq \g$ be a vector space complement of $\g$ and $U_\fe := U \cap \fe$. 
For each $x \in U_\fe$ and $y \in \fn$, we then have $x + y \in U$, so that 
$$ \Spec(\ad_\g(x + y)) \cap 2\pi i \Z \subeq \{0\}. $$
For $x \in \g$, the Lie algebra $\g_x := \fn + \R x$ 
is closed and of finite codimension, 
so that for each $z \in \g_x$ we have the relation $\Spec(\ad_{\g_x} z) 
\subeq \Spec(\ad_\g z)$. Since $\g_x$ is prosolvable, for each root 
$\Gamma(\alpha) \: \Gamma(\g_x) \to \C$ we have 
\begin{eqnarray*}
 \Gamma(\alpha)(\exp_G y \exp_G x) 
&=& \Gamma(\alpha)(\exp_G y)\Gamma(\alpha)(\exp_G x) = \alpha(x)+  \alpha(y) \\
&=& \alpha(x+y) \in \C \setminus (2\pi i \Z\setminus \{0\}).
\end{eqnarray*}
Now Proposition~\ref{prop:4.7} implies that 
$$ \exp_G y \exp_G x \in \exp_G(\fn + x). $$
We conclude that 
$$ \exp_G(\fn) \exp_G(U_\fe) \subeq \exp_G(\fn + U_\fe). $$
Since $\fn \trile \g$ is an ideal, we also have for each $x \in \g$ the relation 
$$ \exp_G(x + \fn) \subeq \exp_G x \exp_G(\fn), $$
which leads to 
$\exp_G(\fn + U_\fe) \subeq \exp_G(\fn) \exp_G(U_\fe),$
and therefore to the equality 
$$\exp_G(\fn + U_\fe) = \exp_G(\fn) \exp_G(U_\fe).$$

The quotient map $\g \to \g/\fn$ integrates to a morphism of Lie groups 
$$ q \: G = \Gamma(\g) \to \Gamma(\g/\fn) $$
whose kernel is  isomorphic to $\Gamma(\fn)$ (\cite{HoMo06}, Theorem~6.7), 
hence equal to $\exp_G(\fn)$ because the exponential map of $\Gamma(\fn)$ 
is surjective. 

We consider the open subset 
$V := \fn + \frac{1}{2} U_{\fe}$
of $\g$. From the preceding considerations it follows that 
$$ \exp(V) = \exp(\fn) \exp(\frac{1}{2} U_\fe). $$
In view of $\ker q = \exp_G(\fn)$, this is the inverse image of an open 
identity neighborhood in $\Gamma(\g/\fn)$, hence an open identity neighborhood in 
$G$. 

For $x,y \in V$, we further know that $\exp_G x = \exp_G y$ implies 
$[x,y] =0$ and $\exp_G (x-y) = \1$. Since $x - y \in \fn + U_{\fe}$ is regular, 
we further get $x - y \in \z(\g)$, and since $\exp_G$ is injective on 
$\z(\g)$, we see that $\exp_G\res_V$ is injective. 
Therefore we have an inverse map $\Psi := (\exp\res_V)^{-1} \: \exp(V) \to V$. 

We want to show that $\Psi$ is smooth. To this end, we write $\g$ 
as a projective limit $\g = \prolim \g_j$ with 
$\g_j \cong \g/\fn_j$, where $\fn_j$ are closed ideals contained in $\fn$. 
Let $\Psi_j := q_j \circ \Psi \: \exp(V) \to \g_j$, where 
$q_j \: \g \to \g/\fn_j$ is the quotient map. Then it suffices to 
show that all maps $\Psi_j$ are smooth. 

The image $q_j(U)$ consists also of regular elements of $\g_j$ 
and $V_j := q_j(V)$ satisfies $V_j + V_j \subeq q_j(U)$, so that 
$$\exp_{\Gamma(\g_j)} \: V_j \to \exp_{\Gamma(\g_j)}(V_j) $$
is a diffeomorphism onto an open subset of $\Gamma(\g_j)$ and 
$\Psi_j$ factors through the inverse map 
$$ (\exp_{\Gamma(\g_j)}\res_{V_j})^{-1} \: 
\exp_{\Gamma(\g_j)}(V_j) = \Gamma(q_j)(\exp(V)) \to V_j. $$
Since $\Gamma(q_j)$ is smooth (Proposition~\ref{prop:cont-smooth}), we 
conclude that 
$$\Psi_j = (\exp_{\Gamma(\g_j)}\res_{V_j})^{-1} \circ \Gamma(q_j) $$
is smooth and hence that the 
Lie group $\Gamma(\g)$ is locally exponential.

(3) $\Rarrow$ (2) has already been observed in Remark~\ref{rem:4.3}(a). 
\end{proof} 

In view of Theorem~\ref{thm:10} we now have 
\begin{corollary} \label{cor:locexp} For a pro-Lie group $G$ the following are equivalent: 
  \begin{description}
  \item[{\rm(1)}] $G$ is locally exponential. 
  \item[{\rm(2)}] $G$ is locally contractible and $\tilde G_0 \cong \Gamma(\fL(G))$ is locally 
exponential.  
  \item[{\rm(3)}]  $G$ is locally contractible and $\fL(G)$ is locally exponential. 
  \end{description}
\end{corollary}

\begin{remark} Let $\g$ be a prosolvable Lie algebra of the form 
$\g = \fn \rtimes_D \R$, where $D \in \der(\fn)$. We assume that 
$\fn$ is exponential, which is equivalent to 
$$ \alpha(\fn) \cap i \R = \{0\}$$ 
for each root $\alpha$ of $\g$. For each nonzero root $\alpha$ of 
$\g$, we put $L_\alpha := \alpha(\fn)$, which is either $\{0\}$ or a 
one-dimensional real subspace of $\C$ intersecting $i\R$ trivially. 

The Lie algebra $\g$ is exponential if and only if we have for each 
root $\alpha$ the stronger condition 
$\alpha(\g) \cap i\R = \{0\}$ (Theorem~\ref{thm:4.1}), which is equivalent to 
$\alpha(D) \in L_\alpha$ for $L_\alpha \not=\{0\}$ and  
$\alpha(D) \not\in i \R$ for $L_\alpha =\{0\}$. 

The condition that all elements in $\fn + [-1,1]D$ are regular is much 
weaker, it means that for each root $\alpha$, we have 
$$ 2\pi i \not \in L_\alpha + [-1,1]\alpha(D). $$
This condition is satisfied in particular if 
$$ |\alpha(D)| < \dist(L_\alpha, 2 \pi i). $$
\end{remark}

\begin{example} Our characterization of the locally exponential 
Lie algebras seems to suggest that 
a pro-Lie algebra $\g$ is locally exponential if and only 
if its radical $\fr = \rad(\g)$ is locally exponential and the representation of 
$\s$ on $\fr$ is ``bounded'' in the sense that its dual contains only finitely 
many types of simple modules. 

The following example shows that this is not the case. 
For each $n \in \N$ we consider the $5$-dimensional real Lie algebra 
$$ \fr_n := \C^2 \rtimes_{D_n} \R, 
\quad \mbox{ where} \quad D_n 
= \begin{pmatrix} 1+ n i & 0 \cr 0 & 1 + ni\end{pmatrix} \in \gl_2(\C). $$
Then each $\fr_n$ is exponential, so that 
$\fr := \prod_{n = 1}^\infty \fr_n$
is an exponential solvable pro-Lie algebra. 

Next we put $\s := \sL_2(\C)$ and let it act on each factor $\Cns^2$ in the 
canonical fashion, which leads to
the semidirect product pro-Lie algebra 
$\g := \fr \rtimes \s,$
whose radical is $\fr$ and for which the representation of $\s$ on $\fr$ is ``bounded''. 

For 
$$ x_n := \Big( \frac{2\pi}{n} D_n,-\frac{2\pi}{n}\begin{pmatrix}1 & 0 \\ 0 & -1\end{pmatrix}\Big) $$
one eigenvalue of $\ad x_n$ on the $n$-th factor space $\Cns^2$ is 
$$  \frac{2\pi}{n}\big((1+ in)-1\big) = 2\pi i, $$
so that $x_n$ is not $\exp$-regular. But in $\g$ we have $x_n \to 0$, 
so that $\g$ is not locally exponential. 
\end{example}

\section{Additional remarks} 

\begin{remark} Let $\g$ be a smooth pro-Lie algebra. 
We have seen above that 
$G := \Gamma(\g)$ is a regular Lie group, hence has a smooth exponential function. 
In view of Example~\ref{ex:4.3}, we cannot expect the group $G$ to 
be locally exponential, that is,  of the first kind in the sense of Robart 
(\cite{Rob97}). 

But we have a class of groups which are still well-behaved and which 
Robart \cite{Rob97} calls {\it groups of the second kind}. Indeed, he defines 
a group to be of the second kind if there exist 
two closed subspaces $\fa, \fb \subeq \L(G)$ such that  
the map 
$$ \fa \times \fb \to G, \quad (x,y) \mapsto \exp_G(x)\exp_G(y) $$
is a local diffeomorphism in $(0,0)$. Considering the differential in this point, 
this requires in particular that $\fa \oplus \fb \cong \L(G)$ as topological 
vector spaces. 

Let us assume that $\g$ contains a closed ideal $\fn$ 
of finite codimension which is exponential and pick a complementary 
subspace $\fe \subeq \g$ (cf.\ Proposition~\ref{prop:4.2}). We claim that the smooth map 
$$ \Phi \: \fn \times \fe \to G, \quad (x,y) \mapsto \exp_G x \exp_G y $$
is a local diffeomorphism in $(0,0)$, hence that $G$ is a Lie group of the second 
kind, regardless of whether it is locally exponential or not. 

Let $\fq := \g/\fn$ be the finite-dimensional quotient Lie algebra, 
$q \: \g \to \fq$ the quotient map, and $\sigma \: \fq \to \fe$ a linear isomorphism. 
Further, let $U \subeq \fq$ be an open $0$-neighborhood such that 
$$\phi := \exp_{\Gamma(\fq)}\res_U \: U \to V := \exp_{\Gamma(\fq)}(U) $$
is a diffeomorphism.

The corresponding homomorphism of groups 
$$ \Gamma(q) \: G = \Gamma(\g) \to Q := \Gamma(\fq) $$
is smooth (Proposition~\ref{prop:cont-smooth}) and 
Theorem~6.7 in \cite{HoMo06} implies that 
$\ker \Gamma(q)$ coincides with the image of the natural map 
$\Gamma(\fn) \to \Gamma(\g),$
which is injective. In view of the exponentiality of $\fn$, 
it follows in particular that 
$$ N := \ker \Gamma(q) = \exp_G(\fn) $$
and that $\exp_G \res_{\fn} \: \fn \to N$ is a bijective smooth map. 
Next we use Gl\"ockner's Implicit Function Theorem (\cite{Gl03}) to see that 
$\ker \Gamma(q)$ is a smooth submanifold of $G$, hence a Lie group 
whose regularity follows from the regularity of $G$ 
(Theorem~\ref{thm:reg-ext}) with the Lie algebra $\ker q = \fn$ (cf.\ \cite{KM97}, 38.7).
Corollary~\ref{cor:reg-cover} now implies that $N \cong \Gamma(\fn)$ as regular 
Lie groups and hence that $\exp_G\res_{\fn} \: \fn \to N$ is a diffeomorphism. 
Let $\log_N := (\exp_G\res_\fn)^{-1} \: N \to \fn$ denote its inverse. 

For $(x,y) \in \fn \times \sigma(U)$, we then have 
$$ \Phi(x,y) = \exp_G(x) \exp_G(\sigma(q(y)), $$
and 
$$ \Gamma(q)(\Phi(x,y)) = \exp_{\Gamma(\fq)}(q(y)) = \phi(q(y)). $$
Hence 
$$ y = \sigma(q(y)) = \sigma(\phi^{-1}(\Gamma(q)(\Phi(x,y)))) $$
implies that $\Phi\res_{\fn \times U}$ is invertible with the inverse 
$$ \Phi^{-1}(g) = (\log_N(g \exp(-\sigma \circ \phi^{-1} \circ \Gamma(q)(g))), 
\sigma \circ \phi^{-1} \circ \Gamma(q)(g)), $$
and this implies that $\Phi\res_{\fn \times U}$ is a diffeomorphism. 
\end{remark}

\begin{remark} \label{rem:pro-smooth} (Regular pro-Lie groups as ``pro-manifolds'') 

Let $G = \Gamma(\g)$ be the regular Lie group associated to 
the smooth pro-Lie algebra $\g$. 

We write $\g = \prolim \g_j$ as a projective limit of finite-dimensional 
Lie algebras $\g_j$ for which the maps $q_j \: \g \to \g_j$ are surjective. 
We have seen in Proposition~\ref{prop:cont-smooth} that 
the corresponding maps 
$\Gamma(q_j) \: G \to G_j := \Gamma(\g_j)$ are smooth. From 
that it follows in particular, that for each map $f \: M \to G$, $M$ a smooth 
manifold, the smoothness of $f$ implies the smoothness of the maps 
$f_j := \Gamma(q_j) \circ f$. 

We claim that the converse also holds, that is,  that $G \cong \prolim G_j$ 
also holds in the category of smooth manifolds. To this end, we have 
to show that the smoothness of all maps $f_j$ implies the smoothness 
of $f$. In Section 3 we have constructed the Lie group structure 
on $G$ by writing it as $R \rtimes S$ and $R \cong \fn \times \fe$, 
hence as a product of three smooth manifolds. Accordingly we write 
$f$ as $f = (\exp_G\circ f_\fn)\cdot  (\exp_G\circ f_\fe) \cdot f_S$, 
and it suffices to show that $f_\fn$, $f_\fe$ and $f_S$ are smooth 
of the corresponding maps into the finite-dimensional quotients are smooth. 

For $f_S$, this follows from the decomposition of $S$ as 
$\tilde\SL_2(\R)^{J_0} \times S_1,$ 
where $S_1$ is finite-dimensional. 
The assertion is obviously true if $\g$ is abelian, which takes care 
of the map $f_\fe$, and, by a straightforward 
projective limit argument, also of $f_\fn$. 
\end{remark}

\begin{remark} (More automatic smoothness) 
An important consequence of Remark~\ref{rem:pro-smooth} is that if 
$H$ is a Lie group and $f \: H \to G = \Gamma(\g)$ is a homomorphism of 
Lie groups, then the smoothness of $f$ follows if the corresponding homomorphisms 
$f_j \: H \to G_j$ are smooth. If, for instance, 
$H$ is locally exponential, 
the continuity of $f$ implies the continuity of all $f_j$, 
and the local exponentiality of the finite-dimensional groups $G_j$ 
implies  that $f_j$ is smooth (\cite{GN06}), hence that $f$ is smooth. 
\end{remark}

\section{Problems} 

\begin{problem} Under which assumptions does continuity of a morphism of pro-Lie groups  
which are Lie groups imply smoothness? Does the existence of a 
smooth exponential function suffice? 

We have the manifold decomposition 
$\Gamma(\g) \cong \fr \times \fe \times \Gamma(\s)$ 
and further $\Gamma(\s) \cong S_0 \times S_1$ with 
$S_1$ finite-dimensional and $S_0 \cong \tilde\SL_2(\R)^{J_0}$. 
In view of  the Iwasawa decomposition $K B$ of $\tilde\SL_2(\R)$, 
we moreover get a diffeomorphism 
$$ S_0 \cong \fk_{\s,0} \times \fb_0, $$
which leads to 
$$ \Gamma(\g) \cong \fr \times \fe \times \fk_{\fs,0} \times \fb_0 \times S_1 $$
as smooth manifolds. 

If $f \: \Gamma(\g) \to H$ is a continuous homomorphism into a Lie group 
$H$ and $H$ has a smooth exponential function, then smoothness of 
$f$ will follow as soon as we have it on each of the five factors above. 
For this it suffices to have a continuous linear map $\phi \: \g \to \h$ 
with $f \circ \exp_{\Gamma(\g)} = \exp_H \circ \phi$, but the existence 
of this map is not obvious in this context. 
\end{problem}

\begin{problem} Let $G$ be a Lie group whose Lie algebra $\L(G)$ is a pro-Lie 
algebra. Show that $G$ has a (smooth) exponential function. 

We know already that for each $x \in \g := \L(G)$ we have a 
smooth $\R$-action on $\g$ generated by the derivation $\ad x$. 
Therefore the main point is to lift this one-parameter group of $\Aut(\g)$ 
through the adjoint action $\Ad \: G \to \Aut(\g)$. 

An important special case is $\g = \R^\N$ and $G$ abelian; but even in 
this case it is not clear how to attack the problem. 
\end{problem}

\section{Appendix. Local contractibility} 

\begin{definition}  \label{def:a.1} (i)  A topological space $X$ is called {\it contractible}
if the identity map of $X$ is homotopic to some constant selfmap of $X$.  
A subspace $Y$ of $X$ is said to be contractible in $X$ to  a
point $y\in Y$ if the inclusion map $Y\to X$ is homotopic 
to the constant map $Y\to X$ with value $y$.

(ii) A topological space $X$ is said to be {\it locally
contractible} at $x\in X$ if there is a neighborhood $U$ 
of $x$ such that $U$ is contractible to $x$ in $X$.

(iii) A homogeneous space $X$ is {\it locally contractible} if it is 
locally contractible at one, and hence at any point.         
\end{definition}

A contractible space is aspherical (that is, has trivial homotopy in all
dimensions) and is acyclic (that is, 
has trivial homology and cohomology with respect to all
homology or cohomology theories that satisfy the Homotopy Axiom).
In particular, $X$ is arcwise connected, and is contractible in $X$ to every 
point $x\in X$.
If $z\in Z\subseteq Y\subseteq X$ and $Y$ is contractible to $z$ in $X$
then $Z$ is contractible to $z$ in $X$.

A space is ostensibly locally contractible at a point $x$ if $x$
has a contractible  neighborhood. 

\begin{remark} \label{rem:a.2} 
Every convex  subset of a locally convex topological vector space
is contractible, and so every open subset of a locally convex topological
vector space is locally contractible. Every manifold modelled
on a locally convex space is locally contractible.
\end{remark}

\begin{lemma}
  \label{lem:a.3} 
Let $X$ and $Y$ be spaces. If $X\times Y$ is locally contractible
at $(x_0,y_0)$, then $Y$ is locally contractible in $y_0$.
If $X$ is contractible  and $Y$ is locally contractible at $y_0$, 
then $X\times Y$ is locally contractible at $(x_0,y_0)$ for
all $x_0\in X$.
\end{lemma}

\begin{proof} First assume that $W$ is a neighborhood of $(x_0,y_0)$ in
$X\times Y$ that is contractible in $X\times Y$. Let 
 let $g_t\colon W\to X\times Y$,
$t\in[0,1]$ be a homotopy such that $g_0(x,y)=(x,y)$ and 
$g_1(x,y)=(x_0,y_0)$ for all $(x,y)\in W$. The set 
$W_Y=\{y\in Y: (x_0,y)\in W\}$ is a neighborhood of $y_0$ in $Y$; let
$h_t\colon W_Y\to Y$ be defined by $h_t(w)=\proj_Y(g_t(x_0,w))$.
Then $(w,t)\mapsto h_t(w): W \times [0,1] \to Y$ is continuous and $h_0(w)=w$
while $h_1(w)=y_0$. Thus $W_Y$ is contractible in $Y$ to $y_0$. 

Now assume that $X$ is contractible. Let $U$ be a neighborhood
 of $y_0$ in $Y$ that is contractible in $Y$ to $y_0$.
Then $X\times U$ is a neighborhood of $(x_0,y_0)$ in $X\times Y$
for each $x\in X$ which is contractible in $X\times Y$.
Indeed, let
$f_s\colon X\to X$ and $g_s\colon U\to Y$, $s\in[0,1]$,
be  homotopies such that $f_0(x)=x$ and
$f_1(x)=x_0$ for all $x\in X$, moreover $g_0(u)=u$ and $g_1(u)=y_0$
for all $u\in U$. Then $(x,u)\mapsto (f_s(x),g_s(u)):X\times U\to X\times Y$,
$s\in[0,1]$, is the required contraction of $X\times U$ to $(x_0,y_0)$
in $X\times Y$.  
\end{proof}


\begin{thebibliography}{aaaaaaaa} 


\bibitem[Bou89]{Bou89} Bourbaki, N., ``Lie Groups and Lie Algebras 
(Chapters 1--3)'', Springer-Verlag, Berlin, 1989. 

\bibitem[Dix57]{Dix57} Dixmier, J., 
{\it L'application exponentielle dans les groupes de Lie r\'esolubles}, 
Bull. Soc. Math. Fr. {\bf 85} (1957), 113-121. 


\bibitem[Gl03]{Gl03} Gl\"ockner, H., {\it Implicit 
functions from topological vector spaces to 
Banach spaces}, Israel J. Math., to appear; math.GM/0303320 

\bibitem[GN06]{GN06} Gl\"ockner, H., and K.-H. Neeb, ``Infinite-Dimensional Lie 
 Groups," book in preparation. 

\bibitem[Hel78]{Hel78} 
Helgason, S., ``Differential Geometry, Lie Groups, and
Symmetric Spa\-ces,'' Acad. Press, London, 1978. 
%
\bibitem[HHL89]{HHL89} Hilgert, J., K.~H. Hofmann, and J.~D.~Lawson, 
``Lie Groups, Convex Cones and Semigroups,'' Oxford, Clarendon Press, 
1989. 

\bibitem[Ho65]{Ho65} Hochschild, G., ``The Structure of Lie Groups,'' Holden Day, San 
Francisco, 1965. 

\bibitem[HoMs66]{HoMs66} Hofmann, K. H., and P. S. Mostert,
``Elements of Compact Semigroups," Charles E. Merrill Books,
Columbus, Ohio, 1966. 

\bibitem[HoMo98]{HoMo98} Hofmann, K.\ H., and S.\ A.\ Morris, ``The Structure of
Compact Groups,'' Studies in Math., de Gruyter, Berlin, 1998. 





\bibitem[HoMo06]{HoMo06} ---, ``{\it The Lie Theory of Connected Pro-Lie Groups--A 
Structure Theory for Pro-Lie Algebras, Pro-Lie Groups and Connected 
Locally Compact Groups}, EMS Publishing House, Z\"urich 2006, to appear.

\bibitem[Iwa49]{Iwa49} Iwasawa, K., {\it On some types of topological groups}, 
Ann. of Math. {\bf 50} (1949), 507--558.

\bibitem[KM97]{KM97} Kriegl, A., and P.\ Michor, ``The Convenient Setting of
Global Analysis,'' Math.\ Surveys and Monographs {\bf 53}, Amer.\
Math.\ Soc., 1997. 

\bibitem[Kur59]{Kur59} Kuranishi, M., {\it On the 
local theory of continuous infinite pseudo groups I}, 
Nagoya Math. J. {\bf 15} (1959), 225--260. 

\bibitem[Lew39]{Lew39} Lewis, D., {\it Formal power series transformations}, Duke Math. J. 
{\bf 5} (1939), 794--805.  

%
%
\bibitem[Mil84]{Mil84} Milnor, J., 
{\it Remarks on infinite-dimensional Lie groups}, pp.~1007--1057,  
In:  B.~DeWitt, R. Stora (eds), 
``Relativit\'{e}, groupes et topologie II 
(Les Houches, 1983), North Holland, Amsterdam, 1984. 

\bibitem[Ne99]{Ne99} Neeb, K.-H., ``Holomorphy and Convexity in Lie Theory," 
Expositions in Mathematics {\bf 28}, de Gruyter Verlag, 1999. 

\bibitem[Ne06]{Ne06} ---, {\em Towards a Lie theory of locally convex 
groups}, Jap. J. Math., to appear. 

\bibitem[Omo80]{Omo80} Omori, H., {\it A method of classifying expansive singularities}, 
J. Diff. Geom. {\bf 15} (1980), 493--512.  

\bibitem[Rob97]{Rob97} Robart, T., {\it Sur l'int\'{e}grabilit\'{e}
des sous-alg\`{e}bres de Lie en dimension infinie}, 
Canad.\ J. Math. {\bf 49:4} (1997), 820--839. 

\bibitem[Sai57]{Sai57} Saito, M., {\it Sur certains groupes de Lie resolubles}, 
Sci. Pap. Coll. Gen. Educ. Univ. Tokyo {\bf 7} (1957), 1--11. 

\bibitem[St61]{Ste61} Sternberg, S., {\it Infinite Lie groups and the formal aspects 
of dynamical systems}, J. Math. Mech. {\bf 10} (1961), 451--474. 


\end{thebibliography}
\end{document}